\documentclass[a4paper,10pt,reqno,oneside]{amsart}    % Europa

\usepackage {graphics}
\usepackage {amsmath}
\usepackage{amsthm}
\usepackage{amssymb}
\usepackage{mathbbol}
\usepackage{color}
\usepackage{url}
\usepackage{enumerate}
\usepackage[round,comma]{natbib}
\usepackage{ifpdf}

\let\cal=\mathcal  %% AMS classes do not define \cal
\let\ve=\relax  %% mkoeppe always types \ve for vectors; make it a no-op

\newtheorem{Theorem}{Theorem}[section]

\newtheorem{Lemma}[Theorem]{Lemma}

\newtheorem{Corollary}[Theorem]{Corollary}
\newtheorem{Definition}[Theorem]{Definition}
\newtheorem{Remark}[Theorem]{Remark}

\newtheorem{Non-Lemma}[Theorem]{Non-Lemma}

\newtheoremstyle{example}{3pt}{3pt}{\normalfont\small}{}{\bfseries\small}{.}{.5em}{}
\theoremstyle{example}
\newtheorem{Example}[Theorem]{Example}

\newcommand {\st} {\mathrm{s.t.}}

\newcommand {\Id} {\mathrm{Id}}
\newcommand {\agg}{\mathrm{agg}}
\DeclareMathOperator\tr{tr}
\DeclareMathOperator\pos{pos}

\newcommand{\R}{{\mathbb R}}

\newcommand{\Z}{{\mathbb Z}}

\newcommand\dotcup{\mathbin{\dot\cup}}

\usepackage%%[breaklinks,colorlinks,citecolor=blue]
{hyperref}

\begin{document}

\title[Integer Minkowski Programs]{Integer Minkowski Programs \\ 
  and the Design of Survivable Networks}

\author{Elke Eisenschmidt}
\address{Elke Eisenschmidt, Otto-von-Guericke-Universit\"at Magdeburg, Department of
  Mathematics, Institute for Mathematical Optimization (IMO),
  Universit\"atsplatz 2, 
  39106 Magdeburg, Germany}
\email{eisensch@imo.math.uni-magdeburg.de}

\author{Matthias K\"oppe}
\address{Matthias K\"oppe, Otto-von-Guericke-Universit\"at Magdeburg, Department of
  Mathematics, Institute for Mathematical Optimization (IMO),
  Universit\"atsplatz 2, 
  39106 Magdeburg, Germany} 
\email{mkoeppe@imo.math.uni-magdeburg.de}
\thanks{The first author was supported by a France T\'el\'ecom research
  agreement with Robert Weismantel.
  The second author acknowledges the support by a Feodor Lynen Research
  Fellowship from the Alexander von Humboldt Foundation and also
  the hospitality of Jes\'us De Loera and the Department of Mathematics of the
  University of California, Davis, where a part of this work was completed.
  The authors wish to thank Robert Weismantel and Raymond Hemmecke for helpful
  discussions. 
}

\author{Alexandre Laugier}
\address{Alexandre Laugier, France T\'el\'ecom R\&D/BIZZ/DIAM, 905 rue Albert Einstein, 06921 Sophia Antipolis
  Cedex, France}
\email{alexandre.laugier@orange-ft.com}

\date{$\relax$Revision: 1.43 $ - \ $Date: 2006/10/27 11:28:47 $ $}

%% MSC2000:
%%   90B10 Network models, deterministic
%% **90B18 Communication networks
%% **68M10 Network design and communication 
%%   90B25 Reliability, availability, maintenance, inspection
%% **90C10 Integer programming
%% **90C11 Mixed integer programming
%% **90C35 Programming involving graphs or networks

\begin{abstract}
  We introduce a new class of optimization problems called \emph{integer
    Minkowski programs}.  The formulation of such problems involves finitely
  many integer variables and nonlinear constraints involving functionals
  defined on families of discrete or polyhedral sets.  We show that, under
  certain assumptions, it is possible to reformulate them as integer linear
  programs, by making use of integral generating sets.  We then apply this
  technique to the network design problem for fractional and integral flows 
  subject to survivability constraints.
\end{abstract}

\maketitle

\section{Introduction}

Let us consider the problem of optimal design of a telecommunications network.
In such a network, digital data are to be sent between pairs of \emph{terminal
  nodes}.  The communicating terminals are, of course, not connected directly
to each other but only indirectly via a network of \emph{switching nodes}.  In
operating and maintaining a telecommunications network, several optimization
problems need to be solved, each on a different time scale.

For the \emph{routing problem}, the network is considered fixed, and certain
\emph{demands} are given.  The direct connections (\emph{links}) between the
switching nodes of a given network have certain given transmission capacities,
usually measured in units of GiBit/s.  The objective of the routing problem is
to find a way to send the data along paths of the network, so that the
transmission capacities on the individual links are not exceeded.  This makes
the routing problem a classical problem of combinatorial optimization, namely
a multicommodity flow problem.  However, the routing problem usually is
complicated by further requirements, like quality-of-service (QoS) guarantees.
For instance, for voice and multimedia communications, it is essential to
bound the \emph{delays}, i.e., the amount of time between
sending and receiving a data packet, along the paths of a routing for the data
of a terminal pair. In these situations, it is often also requested to send
the data of a terminal pair along a unique path of the network, i.e., the flow
cannot be split.

%% Moreover, usually only continuous multicommodity flow problems can be solved
%% exactly.  Here: discrete + hierarchical...

On a larger time scale, of course, the capacities of the links are not fixed.
The capacities are determined both by the characteristics of the 
cables and by the switching nodes. The cables are replaced very rarely because
this is extremely expensive; however the capacities of most links can be
increased by replacing the transmitter--receiver technology installed in the
switching nodes.  We note that, in fiber-optic technology, the installed
capacities are usually symmetric. 
In the \emph{capacitated network design problem}, or more precisely the
\emph{capacity expansion problem}, a set of 
(projected) demands are given, which the current network cannot
satisfy; the objective is to find a set of links and new capacities, such that
the demands can be satisfied, and such that the installation costs are
minimal.  On an even larger time scale, also the addition of new links, changing
the topology of the network, can be planned and implemented, to increase the
capability of the network.  Network design problems have been considered in a
large body of literature, see for instance
\cite{Bienstock,bienstock-chopra-gunluk-tsai-1998,atamturk-2002:network-design-cutset}
and the references within.\bigbreak  

In network design problems, one is often concerned with \emph{survivability}
requirements.  When one link or several links fail, there should still be a
routing that satisfies all demands; this requirement is usually called
\emph{arc survivability}.  Similarly, when a node fails, there should still be
a routing that satisfies all demands.  However, the demands where the failing
node is the source or the destination need not be considered.  This
requirement is called \emph{node survivability}.
Network design problems under node and arc survivability constraints have been
studied in the literature, see for instance
\cite{groetschel-monma-1990,groetschel-monma-stoer-1992,groetschel-monma-stoer-1992b,stoer-1992,stoer-dahl-1994,stoer-dahl-1998,lisser-1995:optimal-joint-synthesis,lisser-1995:survivability,Bienstock2}. 

In the following we will discuss some of the models and solution approaches 
regarding network design problems under survivability constraints. 
To this end, let us fix some notation.
Let $\ve c$ denote a vector of installed capacities, and let $\ve d$ denote a
vector of demands; we shall make the spaces where these vectors live
precise in the following sections.  We shall denote by $\mathcal F_{\ve c,\ve
  d}$ the set of feasible routings.  Again, the precise space that $\mathcal
F_{\ve c,\ve d}$ is a subset of will depend on the choice of the formulation;
we will discuss it later.

\citet{stoer-dahl-1994} seem to have been the first to address survivability
questions in \emph{capacitated} networks; before that paper, research had
focused on modeling graph-theoretic connectivity constraints.  They based
their model on the characterization of capacities that allow
a \emph{fractional routing} of the demands by so-called \emph{metric
  inequalities}
\citep{iri-1971:multi-maxflow-mincut,onaga-kakusho-1971:multi-maxflow-mincut}.
Therefore their model does not contain variables for the no-fault routing and
the fault routings; such a model is known as a \emph{projected formulation}.

There are two shortcomings in this approach.  First, allowing fractional
routings often is too weak a relaxation of the problem.  For the existence of an
integer multicommodity flow, however, there is no known polyhedral description
in general.  Therefore, \emph{no projected formulation} for the case of
integral routings is known in the literature.
Second, more complicated constraints like those coming from QoS
guarantees cannot be incorporated into the model of this approach. For
instance, bounds on the admissible delays along paths of the routing cannot be
modeled, since there are no path variables
present. \citet{stoer-dahl-1994} and various other authors consider
the model of \emph{complete rerouting}, i.e., in the case of a fault the new routing can be
completely unrelated to the no-fault routing.  
%% This can usually not be
%% implemented in practice \citep[cf.][]{lisser-1995:survivability}. 
   
Other models of survivability require  that in the case of a fault only the interrupted flow needs
to be rerouted.  There are two main models here.  In the model of \emph{local
  rerouting}, the failure of a link~$e=(u,v)$ carrying aggregated flow
$x^{\agg{}}_e$ in the no-fault routing creates a new demand to route $x^{\agg{}}_e$
units of flow from node~$u$ to node~$v$ in an auxiliary network where the
failed link~$e$ is not present.  
Local rerouting is also called \emph{link-based re-establishment} in parts of the
literature \citep{Pioro+Medhi:2004}.
In the model of \emph{global rerouting}, 
first the flow along paths that use the failed link is removed from the
no-fault routing; in this way, some capacities on links along these paths are
released.  Then the removed parts of the demands are routed through the 
auxiliary network.  Clearly global rerouting is more economical than local
rerouting, since the released capacities can be used for the fault routings.
Global rerouting is also called \emph{path-based re-establishment} \citep{Pioro+Medhi:2004}.

% \citet{lisser-1995:survivability} consider the problem of computing capacities
% of a dedicated \emph{reserve network} that allows (fractional) global
% rerouting, when the capacities of the \emph{nominal network} and the no-fault
% routing are given.  \citet{lisser-1995:optimal-joint-synthesis} present a
% two-stage method: They first compute a minimum-cost capacity installation 
% of the nominal and the reserve networks that supports a survivable no-fault
% routing, under the assumption of (fractional) local rerouting.  
% In the second stage, they reduce the
% capacities of the reserve network by allowing (fractional) global rerouting.

%%% Mehr Literatur:
%%In \citet{Bienstock2} habe ich das genaue Modell nicht gefunden.
%% Noch interessant waeren evtl. die Referenzen [1,2,18] aus Bienstock2. 

\bigbreak
\subsection*{Contribution of this paper}
%%% FIXME: Goals of the paper?

In this paper we investigate the structure of the family of
feasible-routing sets $\mathcal F_{\ve c,\ve d}$ when the installed capacities~$c$ and the
required demands~$d$ vary.  In our mathematical model, both the capacities
and the demands can take arbitrary (non-negative) integer values.  Thus the
family of feasible-routing sets is infinite.
The key observation in order to be able to deal with this infinite family is that
\emph{Minkowski sums} of these sets play an important r\^ole.  We will show 
that in general the feasible-routing sets satisfy
\begin{equation}
  \mathcal F_{\ve{c^1} ,\ve{d^1}} + \mathcal F_{\ve {c^2},\ve {d^2}}
  \subseteq \mathcal F_{\ve{c^1+c^2}, \ve{d^1+d^2}};\label{eq:minkowski-subseteq}
\end{equation}
then we consider the pairs $(c^1,d^1)$ and $(c^2,d^2)$ where the inclusion is
in fact an equality.

First, it is possible to address the question whether a capacity
installation~$c$ is feasible for a demand vector~$d$.  Whenever there is a
pair $(c^1,d^1)$ and $(c^2,d^2)$ such that $c=c^1+c^2$ and $d=d^1+d^2$ and
in \eqref{eq:minkowski-subseteq} equality holds, then $c$ is feasible for~$d$ 
if and only if $c^1$ is feasible for $d^1$ and $c^2$ is feasible for~$d^2$.

Second, it is possible to address survivability questions as well.  As a first
example, let $a$ be an arbitrary link (arc) of the network.  We will define a
real-valued functional~$g_a$ on the family of feasible-routing sets
by defining $g_a(\mathcal F_{\ve c,\ve d})$ to be the minimum amount of flow carried on~$a$ in
any routing of~$\mathcal F_{\ve c,\ve d}$.  Then the network with capacity installation~$c$ and
demand~$d$ is arc-survivable (in the model of complete rerouting) if and only
if $g_a(\mathcal F_{\ve c,\ve d})=0$ for all arcs~$a$.  It turns out the functional~$g_a$ is
\emph{Minkowski-additive}.  This implies that, if equality holds
in~\eqref{eq:minkowski-subseteq}, then
$(c,d)$ defines a survivable network if and only if $(c^1,d^1)$ and
$(c^2,d^2)$ do. 

\medbreak

\subsection*{Outline}
The outline of the paper is as follows.  In \autoref{IMPs}, we introduce a
class of mathematical optimization problems called \emph{(mixed-)integer
  Minkowski programs} as our general model; their formulation includes
constraints involving non-linear functionals like $g_a$.  We prove a general
reformulation theorem that asserts that mixed-integer Minkowski programs can
be reformulated as integer linear programs, whenever there exists a suitable finite
generating set.  

In \autoref{section atomic fibers} and \autoref{section mixed-integer}, we
show that in two important settings, related to fractional and integral
routing, there exists a finite generating set for the families of feasible-routing
sets~$\mathcal F_{\ve c,\ve d}$.  In \autoref{section formulations} we make
the formulations of the network design problems precise and discuss the
structure of the generating sets that arise.  This enables us to reformulate
the network design problem as an integer linear program.

In \autoref{section survivability}, we address questions of
survivable network design in the model of complete rerouting both for
fractional and integral flows.  Again the survivability conditions can be
modeled using Minkowski-additive functionals in the framework of integer
Minkowski programs.  The general reformulation theorem, applied to this
formulation, then gives a new formulation of the survivable network design
problem as an integer linear program.
The resulting formulation as an integer linear program does not contain
variables for the no-fault routing and the routings in the individual failure
scenarios.  Hence, the formulation can be regarded as a \emph{projected
  formulation}.  This is particularly relevant in the case of integral
routings, where no projected formulation was previously known.  

A prospective advantage of the proposed formulation of the survivable network
design problem is that the formulation only grows moderately with the number
of failure scenarios.  Indeed, only one linear inequality and no extra
variable is needed for each failure scenario.  The number of variables only
depends on the topology of the network.  Moreover, the reformulation only
needs to be computed once for a given network topology; it can then be used
for solving the survivable network design problem for arbitrary demands. 
On the other hand, the usual
non-projected formulations contain one set of routing variables for every
failure scenario.  

%% Finally, in \autoref{section computations}, we show preliminary computational
%% results for a non-trivial example. 

\section{Integer Minkowski programs}

\label{IMPs}

\subsection{Definitions and first results}
\label{Definitions and first results}
In this section we define (mixed-)integer Minkowski programs and show that under
certain conditions there exist reformulations as integer linear programs.
(Mixed-)integer Minkowski programs consist of integer variables, a linear
objective function and inequality constraints involving a certain class of
nonlinear functions.  We first remind the reader of the definition of
Minkowski sums.
\begin{Definition}~
  \begin{enumerate}[\rm (i)]
  \item Let $A_1,A_2\subseteq\R^d$ be sets.  The \emph{Minkowski sum}
    (pointwise sum) of $A_1$ and $A_2$ is defined as 
    \begin{displaymath}
      A_1+A_2 = \bigl\{\, x^1 + x^2 : x^1\in A_1,\, x^2\in A_2\, \bigr\}.
    \end{displaymath}
  \item Let $A\subseteq\R^d$ be a set and let $k\in\Z_+$ be a
    non-negative integer.  We will denote by $k A$ the set obtained by
    taking the Minkowski sum of $k$ copies of the set~$A$,
    \begin{equation}
      \label{eq:k-fold-minkowski}
      k A = \bigl\{\, x^1+\dots+x^k : x^1,\dots,x^k\in A\,\bigr\}.
    \end{equation}
    In the same way, we will use the notation $\sum_{i\in I} k_i A_i$. 
  \end{enumerate} 
\end{Definition}
Note that when one of the summands is the empty set, then the
Minkowski sum is empty too.

\begin{Remark}
  Frequently, Minkowski sums are considered in the context of convex
    sets only.  For convex sets, the notation $kA$ is usually defined for
    arbitrary scalars $k\in\R$ as 
    \begin{displaymath} 
      \bigl\{\, k x : x\in A\,\bigr\};
    \end{displaymath}
    for $k\in\Z_+$ this gives the same set as in~\eqref{eq:k-fold-minkowski}.
    However, when $A$ is not a convex set, the set $kA$ as defined
    in~\eqref{eq:k-fold-minkowski} is not necessarily the same as the set
    $\bigl\{\, k x : x\in A\,\bigr\}$. 
\end{Remark} 

We will consider a class of nonlinear functions that are obtained in a two-stage
process. 
In the first stage, we define a set-valued mapping ${\cal A}\colon \Z^n
\rightarrow 2^{\R_+^d}$ that is superadditive:
\begin{Definition}
  A set-valued mapping ${\cal A}\colon \Z^n \rightarrow 2^{\R_+^d}$ is called
  superadditive if it satisfies 
  \begin{equation}
    \label{inclusion}
    {\cal A}(z^1 + z^2) \supseteq {\cal A}(z^1) + {\cal A}(z^2) \quad
    \text{for } z_1,z_2 \in \Z^n,
  \end{equation}
  where $+$ in the right-hand side denotes a Minkowski sum.
\end{Definition}
In the second stage, we define a functional~$\kappa$ on the range
of~$\mathcal A$ with the following
properties: 
\begin{Definition}
Let $\kappa\colon 2^{\R_+^d} \supseteq \mathcal A(\Z^n) 
\rightarrow \bar{\R} = \R \cup \{+\infty\}$ be a
functional with $\kappa (\emptyset) = +\infty$ and $A,B,C \in \mathcal A(\Z^n)$.
\begin{enumerate}[\rm (i)]
\item $\kappa$ is called \emph{Minkowski-additive} if $C = A +
  B$ implies $\kappa(C) = \kappa(A) + \kappa(B)$.
\item $\kappa$ is called \emph{antitone} if $A \subseteq B$
  implies $\kappa(B) \leq \kappa(A)$.
\end{enumerate}
\end{Definition}

With these preparations, we are now in the position to define (mixed-)integer
Minkowski programs. 

\begin{Definition}
A \emph{(mixed-)integer Minkowski program} is an optimization problem of the
following structure:
\begin{equation}
\label{IMP}
\begin{aligned}
\min\quad & W^T z \\
\st\quad & {\cal A}(z) \neq \emptyset \\
& B \cdot z = v\\
& \kappa^i({\cal A}(z)) + h_i^Tz \leq \delta_i &&\text{for $i=1, \ldots, k$}\\
& z \in \Z^n,\\
\end{aligned}
\end{equation}
where $W \in \R^n$ is a cost vector,
$\mathcal A\colon\Z^n\to 2^{\R_+^d}$ is a superadditive set-valued
mapping,
$B \in \Z^{n \times n}$ a matrix, $v \in \Z^n$,
$\kappa^i\colon \mathcal A(\Z^n)\to \bar{\R}$ for $i=1, \ldots, k$
is a Minkowski-additive and antitone functional, and $h_i \in \R^n$ for
$i=1,\dots,k$.

In the mixed-integer case, we have ${\mathcal A}(z) \subseteq \R_+^d$
for all $z \in \Z^n$, whereas in the integer case ${\mathcal A}(z)
\subseteq \Z^d_+$ for all $z\in \Z^n$.
\end{Definition} 

The following theorem will concretize the conditions under which we
can reformulate a (mixed-)integer Minkowski program as an integer
linear program.

\begin{Theorem}
\label{IP und IMP}
Let $\{\bar{z}^j\}_{\mathcal J} \subseteq \Z^n$ with ${\mathcal
  A}(\bar{z}^j) \neq \emptyset$ for all $j\in {\mathcal J}$ be a
finite  generating set of the 
family of sets ${\mathcal A}(z)$, i.e., for $z\in \Z^n$ there is a
representation 
\begin{displaymath}
  z = \sum\limits_{j \in {\mathcal J}}{\lambda_j \, \bar{z}^j}
  \quad\text{where}\quad \lambda_j \in \Z_+,
\end{displaymath}
such that
$$ {\mathcal A}(z)={\mathcal A}\biggl(\sum\limits_{j \in {\mathcal
    J}}{\lambda_j \, \bar{z}^j}\biggr) = \sum\limits_{j \in {\mathcal
    J}}{\lambda_j \, {\mathcal A}(\bar{z}^j)}. $$
Then the (mixed-)integer Minkowski program \eqref{IMP} can be reformulated
as an integer linear program:
\begin{equation}
\label{IP}
\begin{aligned}
  \min\quad & \sum\limits_{j \in {\mathcal J}} (W^T \bar{z}^j) \lambda_j\\
  \st\quad  & \sum\limits_{j \in{\mathcal J}}{ (B 
    \bar{z}^j) \lambda_j} = v\\ 
& \sum\limits_{j \in {\cal J}}  \bigl[\kappa^i({\cal A}(\bar{z}^j)) +
  h_i^T\bar{z}^j \bigr] \lambda_j  \leq \delta_i && \text{for $i=1, \ldots, k$}\\
 & \lambda_j \in \Z_+ && \text{for $j \in {\cal J}$}\\
\end{aligned}
\end{equation}
\end{Theorem}
\begin{proof}
Let $\lambda_j^*$, $j \in {\mathcal J}$ be a feasible solution of the 
integer program \eqref{IP}. We need to prove that $z^* =
\sum\nolimits_{j \in {\cal J}}{\lambda_j^*\, \bar{z}^j}$ is a feasible 
solution of the corresponding (mixed-)\penalty0 integer Minkowski program
\eqref{IMP}. Relation \eqref{inclusion} yields: 

\begin{equation}
\label{incl}
{\cal A}(z^*) = {\cal A} \biggl(\sum\limits_{j \in {\cal J}}{\lambda_j^* \,\bar{z}^j}\biggr)
\supseteq \sum\limits_{j \in {\cal J}}{\lambda^*_j \, {\cal
    A}(\bar{z}^j)},
\end{equation}
and as ${\cal A}(\bar{z}^j) \neq \emptyset$ this implies ${\cal A}(z^*)
\neq \emptyset$. We also have $B \cdot z^* = v$. It remains to prove that the linear
constraints of the (mixed-)integer Minkowski program  are
satisfied. Let $i \in \{1,
\ldots, k\}$.  As $\kappa^i$ is an antitone functional for all $i=1,
\ldots, k$ and as \eqref{incl} holds, we have: 
$$\kappa^i({\cal A}(z^*)) + h_i^Tz^* \leq \sum\limits_{j \in {\cal J}}{\lambda^*_j\,
  \kappa^i({\cal A}(\bar{z}^j))} + \sum\limits_{j \in {\cal J}}{\lambda^*_j\,
  h_i^T\bar{z}^j}
 \leq \delta_i . $$
This proves the feasibility of $z^* = \sum\limits_{j \in {\cal J}}{\lambda^*_j \,
  \bar{z}^j}$ for~\eqref{IMP}. 

On the other hand, let $\tilde{z}$ be a feasible solution of the
(mixed-)integer Minkowski program \eqref{IMP}. We know that there is a
decomposition $\tilde{z} = \sum\nolimits_{j \in {\cal
    J}}{\beta_j\,\bar{z}^j}$ such that 
$${\cal A}(\tilde{z}) = \sum\limits_{j \in {\cal J}}{\beta_j \, {\cal
    A}(\bar{z}^j)},$$ 
where $\beta_j \in \Z_+$ for all $j \in {\mathcal J}$.
We shall prove that $\beta_j, \; j \in {\mathcal J}$, is a feasible
solution of the integer linear program \eqref{IP}. Surely we have $B
\cdot \tilde{z} = \sum\nolimits_{j \in{\mathcal J}}{ (B 
    \bar{z}^j) \beta_j} = v$.  Let $i \in \{1,
\ldots, k\}$ be fixed. As
$\kappa^i$ is a Minkowski-additive functional for all $i=1, \ldots, k$ we have
$$\delta_i  \geq  \kappa^i({\cal A}(\tilde{z})) + h_i^T\tilde{z}
 =  \sum\limits_{j \in {\cal J}}{\beta_j\,[\kappa^i({\cal A}(\bar{z}^j)) +
  h_i^T\bar{z}^j ]}.$$
Therefore, $\beta_j, \; j \in {\mathcal J}$ is a feasible solution of
the  integer linear program \eqref{IP}. 

Clearly also the objective value of the integer linear program \eqref{IP} is
equal to the objective value of the corresponding (mixed-)integer Minkowski
program \eqref{IMP} for every cost vector $W \in \R^n$.  \end{proof}

\subsection{Families of truncated sets and truncated integer Minkowski programs}
The results of the previous subsection may be improved in the sense
that it is possible to formulate integer Minkowski programs on
families of truncated sets as integer linear program. This result is
important for our application in network design because there are usually
lots of circulations of flow as feasible solutions for the
multicommodity flow problem which are not used for a ``regular''
objective function and therefore should be cut away. We give the
definition of a truncated set first and continue giving the results
for truncated (mixed-)integer Minkowski programs.

\begin{Definition}
\label{sqsubseteq}
Let $u, v \in \R^n$ be vectors. We say that $u$ \emph{reduces} $v$ and
denote $u \sqsubseteq v$ if $u^{(i)}v^{(i)} \geq 0$  and $|u^{(i)}|
\leq |v^{(i)}|$ for all components $i=1, \ldots, n$.
\end{Definition}

\begin{Definition}
\label{truncated sets}
Let ${\mathcal C} \subseteq \R^n$ (or ${\mathcal C} \subseteq \Z^n$)
be a set of vectors. Let ${\mathcal A} \colon \Z^n \rightarrow
2^{\R^d_+}$ be a set-valued mapping which is super-additive and $z
\in \Z^n$. We call
\begin{displaymath}
\tr_{\mathcal C}({\mathcal A}(z)) = \{ x \in {\mathcal A}(z) :
\nexists c \in {\mathcal C} \text{ with } c \sqsubseteq x\}
\end{displaymath}
the truncated set corresponding to $z$ and ${\mathcal A}$ with
respect to ${\mathcal C}$.
\end{Definition}

\begin{figure}[ht]
\centering
\input{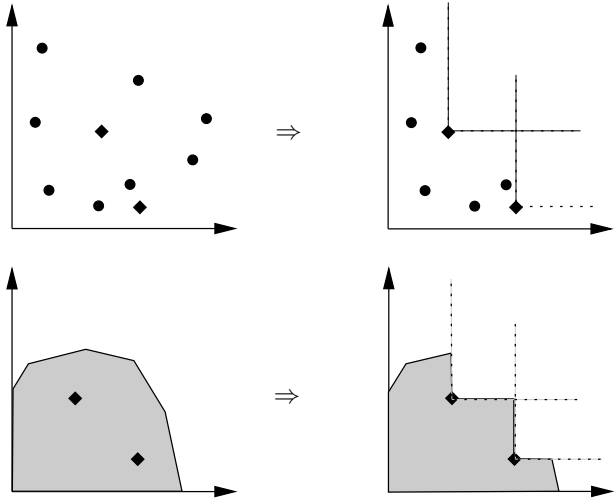}
\caption{Truncation for a discrete set of vectors and a polytope}
\label{truncation-example}
\end{figure}

The definition of truncated sets is illustrated in Figure
\ref{truncation-example} for a discrete set of 
vectors and a polytope. Both are truncated by a fixed set ${\mathcal
  C}$ of two vectors.

\begin{Definition}
\label{indecomposability truncation}
A truncated set $\tr_{\mathcal C}({\mathcal A}(z))$ is called
indecomposable if there is no $0 \neq z_1, z_2 \in \Z^n$ with
\begin{displaymath}
\tr_{\mathcal C}({\mathcal A}(z)) \subseteq \tr_{\mathcal C}({\mathcal
  A}(z_1)) + \tr_{\mathcal C}({\mathcal A}(z_2)).
\end{displaymath}
\end{Definition}

\begin{Lemma}
\label{trunc subset normal}
If $\tr_{\mathcal C}({\mathcal A}(z))$ is indecomposable then
${\mathcal A}(z)$ is indecomposable, too.
\end{Lemma}

\begin{proof}
We have $\tr_{\mathcal C}({\mathcal A}(z)) \subseteq {\mathcal
  A}(z)$. Suppose there is a decomposition ${\mathcal A}(z)=
{\mathcal A}(z_1) + {\mathcal A}(z_2)$. Let $x \in \tr_{\mathcal
  C}({\mathcal A}(z)) \subseteq {\mathcal A}(z)$, then there is a
decomposition $x = x_1 + x_2$ with $x_i \in {\mathcal A}(z_i)$,
$i=1,2$. In particular: $x_i \sqsubseteq x$ (because $x_i \in
\R_+^d$). Suppose there is $c \in {\mathcal C}$ with (w.l.o.g.) $c
\sqsubseteq z_1$, implying $c \sqsubseteq z$, which is a
contradiction. Thus $z_i \in \tr_{\mathcal C}({\mathcal A}(z_i))$ and
we have $\tr_{\mathcal C}({\mathcal A}(z)) \subseteq \tr_{\mathcal
  C}({\mathcal A}(z_1)) + \tr_{\mathcal C}({\mathcal A}(z_2))$,
contradicting the indecomposability of $\tr_{\mathcal C}({\mathcal
  A}(z))$. \end{proof}

The previous lemma shows that there is a finite generating system of
families of truncated sets provided there is such a finite generating
system for the corresponding family of sets. We go on showing that
under certain conditions it is possible to reformulate integer
Minkowski programs on truncated sets as integer linear programs.

Furthermore, let $W \in \R^n$ denote a vector of costs, $B \in \Z^{n
  \times n}$ a matrix, $v \in \Z^n$, 
$\kappa^i\colon \mathcal A(\Z^n)\to \bar{\R}$ for $i=1, \ldots, k$
is a Minkowski-additive and antitone functional, and $h_i \in \R^n$ for
$i=1,\dots,k$.
We consider the following truncated integer Minkowski program:

\begin{equation}
\label{tr-IMP}
\begin{aligned}
\min\quad & W^T z \\
\st\quad & \tr_{\mathcal C}({\mathcal A}(z)) \neq \emptyset \\
& B \cdot z = v\\
& \kappa^i(\tr_{\mathcal C}({\mathcal A}(z))) + h_i^Tz \leq \delta_i &&\text{for $i=1, \ldots, k$}\\
& z \in \Z^n.\\
\end{aligned}
\end{equation}

\begin{Lemma}
\label{tr-reformulation}
Let $\{\bar{z}_j\}_{j \in {\mathcal J}}$ denote the finite generating
system of a family of truncated sets, i.e.,~for $z \in \Z^n$ there are
$\lambda_j \in \Z_+$, $j \in {\mathcal 
  J}$, such that $z = \sum\limits_{j \in {\mathcal J}}{\lambda_j \bar{z}_j}$
and
\begin{displaymath}
\tr_{\mathcal C}({\mathcal A}(z)) \subseteq \sum\limits_{j \in
  {\mathcal J}}{\lambda_j \tr_{\mathcal C}({\mathcal A}(\bar{z}_j))}
\end{displaymath}
If $\kappa^i(\tr_{\mathcal C}({\mathcal A}(z))) = \kappa^i({\mathcal A}(z))$ for
all $z \in \Z^n$ and $i =1, \ldots, k$ then we can reformulate the
truncated integer Minkowski program \eqref{tr-IMP} as integer linear program:
\begin{equation}
\label{tr-IP}
\begin{aligned}
  \min\quad & \sum\limits_{j \in {\mathcal J}} (W^T \bar{z}^j) \lambda_j\\
  \st\quad  & \sum\limits_{j \in{\mathcal J}}{ (B 
    \bar{z}^j) \lambda_j} = v\\ 
& \sum\limits_{j \in {\cal J}}  \bigl[\kappa^i(\tr_{\mathcal
  C}({\mathcal A}(\bar{z}^j))) +
  h_i^T\bar{z}^j \bigr] \lambda_j  \leq \delta_i && \text{for $i=1, \ldots, k$}\\
 & \lambda_j \in \Z_+ && \text{for $j \in {\cal J}$}\\
\end{aligned}
\end{equation}
\end{Lemma}

\begin{proof}
Let $z$ be a feasible solution of the truncated integer Minkowski program
\eqref{tr-IMP}. Then there are $\beta_j \in \Z_+$ such that $z =
\sum\limits_{j \in {\mathcal J}}{\beta_j \, \bar{z}^j}$ and 
$$\tr_{\mathcal C}({\mathcal A}(z)) \subseteq \sum\limits_{j \in {\mathcal
    J}}{\beta_j \, \tr_{\mathcal C}({\mathcal A}(\bar{z}^j))}. $$
It is clear that $\sum\limits_{j \in{\mathcal J}}{ (B 
    \bar{z}^j) \lambda_j} = v$ holds and as the functions $\kappa^i$
  are antitone for $i=1, \ldots, k$ we have:
\begin{displaymath}
\delta_i \geq \kappa^i(\tr_{\mathcal C}({\mathcal A}(z))) + h_i^Tz \geq
\sum\limits_{j \in {\mathcal J}} \bigl[\kappa^i(\tr_{\mathcal
  C}({\mathcal A}(\bar{z}^j))) + h_i^T \bar{z}^j \bigr] \beta_j 
\end{displaymath}
for $i=1, \ldots, k$. This means that $\beta_j \in \Z_+$ for $j \in
{\mathcal J}$ is a feasible solution of the integer linear program
\eqref{tr-IP}. 

It remains to show the other direction. Let $\lambda_j \in \Z_+$ for
$j \in {\mathcal J}$ be a feasible solution of the integer linear
program \eqref{tr-IP}. We set $z:= \sum\limits_{j \in {\mathcal
    J}}{\lambda_j \, \bar{z}^j}$ and show that $z$ is a feasible
solution of the truncated integer Minkowski program \eqref{tr-IMP}. We have:
$$v = \sum\limits_{j \in{\mathcal J}}{ (B 
    \bar{z}^j) \lambda_j} = B \cdot z.$$
On the other hand, as  $z = \sum\limits_{j \in {\mathcal
    J}}{\lambda_j \, \bar{z}^j}$ we have ${\mathcal A}(z) \supseteq
\sum\limits_{j \in {\mathcal J}}{\lambda_j \, {\mathcal A}(\bar{z}^j)}$. This implies:
\begin{equation}
\notag
\begin{split}
\delta_i &\geq \sum\limits_{j \in {\cal J}}
\bigl[\kappa^i(\tr_{\mathcal C}({\mathcal A}(\bar{z}^j))) +
  h_i^T\bar{z}^j \bigr] \lambda_j\\
&= \sum\limits_{j \in {\mathcal J}}
\bigl[\kappa^i({\mathcal A}(\bar{z}^j)) + h_i^T \bar{z}^j \bigr] \lambda_j\\
&\geq \kappa^i({\mathcal A}(z)) + h_i^Tz \\
&= \kappa^i(\tr_{\mathcal C}({\mathcal A}(z))) + h_i^Tz,
\end{split}
\end{equation}
for $i=1, \ldots, k$. This means that $z$ is a feasible solution of
the truncated integer Minkowski program \eqref{tr-IMP}. \end{proof}

In the following sections, we will investigate some finitely generated families of
sets ${\mathcal A}(z)$ and the structures of the
associated truncated (mixed-)integer Minkowski programs.

\section{The integer case: Atomic fibers}
\label{section atomic fibers}
In this section we will investigate the situation when the mapping~$\mathcal
A$ generates a family of \emph{discrete} sets.  In our application, the
discrete sets will be sets of feasible integral routings.  To be more precise,
let $A\in\Z^{d\times n}$ be a matrix, and let the mapping $\mathcal
A\colon\Z^d \to 2^{\R^n_+}$ be 
defined by
\begin{equation}
  \mathcal A(b) = P^I_{A,b} = \{\,z : Az = b,\, z \in \Z^n_+ \,\}.
\end{equation}
The set $P^I_{A,b}$ is known as the \emph{fiber} of~$b$ under the linear map
$f_A\colon \Z^n_+ \rightarrow \Z^d$, $x \mapsto Ax$.

\citet{EHK:2006} considered the family of fibers of a fixed matrix $A \in
\Z^{d \times n}$, when the right-hand side vector $b \in \Z^d$ varies.  
They established a theory of Minkowski decomposition for fibers; in this
context, the non-decomposable fibers are called \emph{atomic}:

\begin{Definition}
  \label{atomic fiber} We call a fiber $P_{A,b}^I$ \emph{atomic}, if
  $P_{A,b}^I = P_{A,b_1}^I + P_{A,b_2}^I$ implies $b = b_1$ or $b =
  b_2$.
\end{Definition}

%% \texttt {We will remove the following subsection from the paper because
%%   all the results are shown in our paper \citet*{EHK:2006}, which
%%   will supersede the paper \citet{Hemmecke}.} 

For our purposes, we will need to consider a slight generalization of the
notion of atomic fibers, where we restrict the set of right-hand sides from
$\Z^d$ to some subset.  As a motivating example, consider flow-conservation
constraints that appear in standard node-arc formulations.  These
constraints are linear equations $A_{i,\cdot} z = b_i$ with a right-hand
side~$b_i=0$.  It is desirable that the flow-conservation constraints also
hold in the Minkowski summands $P^I_{A,b^1}$, $P^I_{A,b^2}$ of a set
$P^I_{A,b}$ of 
feasible integral routings.  Therefore, we wish to restrict the decomposition
into atomic fibers such that $b_i=0$ for all components corresponding to the
flow-conservation constraints.

This example suggests to restrict the right-hand sides to a sublattice
of~$\Z^d$, but later we will see that it is also useful to consider
restriction to a submonoid of~$\Z^d$.
\begin{Definition}
\label{monoid-atomic}
Let $(M,+)$ with $M \subseteq \Z^d$ be a monoid, $(\Lambda, +)$ a
lattice with $\Lambda \subseteq \Z^d$, and let $A \in \Z^{d
  \times n}$ be a matrix.
\begin{enumerate}
\item  A fiber $P_{A,b}^I$ with $b \in \Lambda$ is called
\emph{atomic w.r.t.~$\Lambda$} if there is no decomposition $b= b^1 + b^2$
with $P_{A,b}^I = P_{A,b^1}^I + P_{A,b^2}^I$, where $b^1, b^2 \in \Lambda$.
\item A fiber $P_{A,b}^I$ with $b \in M$ is called
\emph{atomic w.r.t.~$M$} if there is no decomposition $b= b^1 + b^2$
with $P_{A,b}^I = P_{A,b^1}^I + P_{A,b^2}^I$, where $b^1, b^2 \in M$.
\end{enumerate}
\end{Definition}

Once again we refer to the theory of indecomposable fibers
w.r.t.~monoids and 
lattices developed by \citet{EHK:2006} and recall a result which
is important for our analysis in this paper:

\begin{Lemma}
Let $M \subseteq \Z^n$ be a finitely generated monoid and let $A \in
\Z^{d \times n}$ be a matrix. There are only finitely many fibers of
the matrix A that are atomic w.r.t.~$M$.
\end{Lemma}
 
Note that for $M= \langle \pm A_1, \ldots, \pm A_n\rangle$, where
$A_i$, $i=1, \ldots, n$, denote the columns of matrix $A$, the above lemma shows that
there are only finitely many atomic fibers of a matrix $A$.

\subsection{Truncated fibers}

The usual formulations of the routing problem (multicommodity flow problem)
include routings that contain certain flow circulations as feasible solutions.
When a useful (``regular'')
objective function is chosen, an optimal solution to the routing
problem will never contain a flow circulation, when we 
consider the routing commodity by commodity.  The reason is that a routing
with circulations will always be dominated by the routing where the
circulations have been removed.  For our purposes, however, it makes sense to
cut away circulations explicitly; we will see later that it can significantly
simplify the computations.

This application gives rise to the notion of \emph{truncated fibers}.
\begin{Definition}
\label{truncated fibers}
Let ${\mathcal C} = \{c_1, \ldots, c_s\}$ with $c_i \in \Z_+^n$ be a finite
set of vectors.  Let $P_{A,b}^I$ be a fiber of a matrix $A \in \Z^{d \times
  n}$.  We call
$$\tr_{\mathcal C}(P_{A,b}^I) = \{z \in P_{A,b}^I: \; \nexists c
\in {\mathcal C} \text{ with } c \sqsubseteq z\}$$
the \emph{truncated fiber} of $A$ with right-hand side $b$ with respect to
${\mathcal C}$.
\end{Definition}

\begin{Definition}
\label{atomic truncated fibers}
A truncated fiber $\tr_{\mathcal C}(P_{A,b}^I)$ is called
indecomposable (w.r.t.~a monoid $M$) if there is no $b_1, b_2 \neq 0$
($b_1, b_2 \in M$) with $b = b_1 + b_2$ and
$$\tr_{\mathcal C}(P_{A,b}^I) \subseteq \tr_{\mathcal C}(P_{A,b_1}^I) +
\tr_{\mathcal C}(P_{A,b_2}^I).$$
\end{Definition}

\begin{Corollary}
If $\tr_{\mathcal C}(P_{A,b})^I$ is indecomposable (w.r.t.~a monoid $M$), then
$P_{A,b}^I$ is atomic (w.r.t.~the monoid), too.
\end{Corollary}
\begin{proof}
This is a direct consequence of Lemma \ref{trunc subset normal}. \end{proof}

This leads to the main result of this section: We can reformulate the
integer Minkowski program on truncated fibers as an integer linear
program. Consider the integer Minkowski program:
\begin{equation}
\label{tr-IMP-fiber}
\begin{aligned}
\min\quad & W^T~b \\
\st\quad & \tr_{\mathcal C}(P_{A,b}^I) \neq \emptyset \\
& B \cdot b = u\\
& \kappa^i(\tr_{\mathcal C}(P_{A,b}^I) + h_i^Tb \leq \delta_i &&\text{for $i=1, \ldots, k$}\\
& b \in \Z^d,\\
\end{aligned}
\end{equation}
where $W \in \R^d$ is a vector of costs, $B \in \Z^{d \times d}$, $u
\in \Z^d$, $\kappa^i$ antitone and Minkowski-additive functionals,
$h_i \in \R^d$ and $\delta_i \in \R$ for $i=1, \ldots, k$. Then we
have:

\begin{Corollary}
\label{reformulating truncated fibers}
If $\kappa^i(\tr_{\mathcal C}(P_{A,b}^I) = \kappa^i(P_{A,b}^I)$ for
all $b \in \Z^d$ and $i =1, \ldots, k$ then we can reformulate the
integer Minkowski program \eqref{tr-IMP-fiber} as an integer linear program:
\begin{equation}
\label{tr-IP-fiber}
\begin{aligned}
  \min\quad & \sum\limits_{j \in {\mathcal J}} (W^T \bar{b}^j) \lambda_j\\
  \st\quad  & \sum\limits_{j \in{\mathcal J}}{ (B 
    \bar{b}^j) \lambda_j} = u\\ 
& \sum\limits_{j \in {\cal J}}  \bigl[\kappa^i(\tr_{\mathcal
  C}(P_{A,\bar{b}^j}^I) +
  h_i^T\bar{b}^j \bigr] \lambda_j  \leq \delta_i && \text{for $i=1, \ldots, k$}\\
 & \lambda_j \in \Z_+ && \text{for $j \in {\cal J}$}\\
\end{aligned}
\end{equation}
\end{Corollary}
\begin{proof} 
  This is a direct consequence of Lemma \ref{tr-reformulation}. 
\end{proof} 
\subsection{Projections of fibers}
The two basic formulations for the routing problem are the node-arc
formulation and the path formulation.  It turns out that the node-arc
formulation essentially is a projection of the path formulation. By
the term projection we mean linear, integral transformations of problems
from a high dimensional space of variables to a space of variables of
lower dimension.  

Because of this relation between the node-arc and path formulations, we also need
to consider projections of (atomic) fibers.  We consider the fibers of the
integral matrix $A \in \Z^{d \times n}$ under a projection described by a
matrix $\Pi\in \Z_+^{m \times n}$.  We remark that projections of fibers are
not necessarily fibers of some matrix~$\bar A$; they can have a more
complicated structure, see \cite{Williams-1992}.  In our application, however,
the following setting is general enough. 
We suppose that there
is an integral transformation matrix $N \in \Z^{c \times d}$ for the
set of right-hand sides $b \in \Z^d$, with the columns of $N$ being
linearly independent. Furthermore we suppose that
there is an an integral matrix 
$\bar{A} \in \Z^{c \times m}$ such that $N  A = \bar{A} 
\Pi$.  Then we have the following property:
\begin{equation}
\label{identity of fibers}
\begin{split}
\Pi \cdot P_{A,b}^I = \Pi \cdot \{z \in \Z^n_+ : Az=b\} &= \Pi \cdot
\{z \in \Z^n_+ : N  A  z = N  b\}\\
&= \Pi \cdot \{z \in \Z^n_+ : \bar{A}  \Pi z = N
b\}\\
&= \{y \in \Z^m_+ : \bar{A} y = \bar{b}\} = P_{\bar{A},\bar{b}}^I,
\end{split}
\end{equation}
where $\bar{b} = N  b$.
Equation \eqref{identity of fibers} gives the following lemma.

\begin{Lemma}
\label{projected fibers}
If $P_{A,b}^I$ is decomposable then $\Pi \cdot P_{A,b}^I =
P_{\bar{A},\bar{b}}^I$ is decomposable w.r.t.~the monoid ${\mathcal B}
= \{\bar{b} = N \cdot b : b \in A\Z^n_+\}$.
\end{Lemma}
\begin{proof}
This is a direct consequence of the linearity the projection $\Pi$:
\begin{equation}
\notag
\begin{split}
P_{\bar{A},\bar{b}}^I = \Pi \cdot P_{A,b}^I &= \Pi \cdot (P_{A,b_1}^I +
P_{A,b_2}^I)\\
&= \Pi \cdot P_{A,b_1}^I + \Pi \cdot P_{A,b_2}^I\\
&= P_{\bar{A},\bar{b_1}}^I + P_{\bar{A},\bar{b_2}}^I.
\end{split}
\end{equation}
\end{proof}

An analogous assertion is true for truncated fibers.

\begin{Lemma}
\label{truncated,projected fibers}
 Let
${\mathcal C} = \{c_1, \ldots, c_t\} \subseteq \Z^n_+$ a finite set of
vectors and let 
${\mathcal D} \subseteq \Z_+^m$ such that $\bar{A} \cdot d \in
{\mathcal B}$ for all $d \in {\mathcal D}$, where $\mathcal B$ is the monoid
from Lemma~\ref{projected fibers}.  For
all $c \in {\mathcal C}$ let there exist $d \in {\mathcal D}$ with $d
\sqsubseteq \pi(c)$.  Then we have:
If $\tr_{\mathcal C}(P_{A,b}^I)$ is decomposable, then $\tr_{\mathcal
  D}(P_{\bar{A},\bar{b}}^I)$ is decomposable w.r.t.~${\mathcal B}$.
\end{Lemma}
\begin{proof}
As a first step we will show that $\tr_{\mathcal
  D}(P_{\bar{A},\bar{b}}^I) = \tr_{\mathcal D}(\Pi \cdot P_{A,b}^I)
\subseteq \Pi \cdot \tr_{\mathcal C}(P_{A,b}^I)$. To 
this aim let $x \in  \tr_{\mathcal D}(\Pi \cdot P_{A,b}^I)$. Then there is
$y \in P_{A,b}^I$ with $\Pi \cdot y = x$. Suppose there is $c \in
{\mathcal C}$ with $c \sqsubseteq y$. Then there is $d \in {\mathcal
  D}$ with $d \sqsubseteq \Pi \cdot c \sqsubseteq \Pi \cdot y = x$ which
contradicts the fact that $x \in \tr_{\mathcal
  D}(P_{\bar{A},\bar{b}}^I)$. Thus $y \in \tr_{\mathcal
  C}(P_{A,b}^I)$ and consequently $x = \Pi \cdot y \in \Pi \cdot
\tr_{\mathcal C}(P_{A,b}^I)$. This proves our first claim. 
Now let $\tr_{\mathcal C}(P_{A,b}^I)$ be decomposable, i.e., there are
vectors $b_1,b_2$ 
with $\tr_{\mathcal C}(P_{A,b}^I) \subseteq \tr_{\mathcal C}(P_{A,b_1}^I) +
 \tr_{\mathcal C}(P_{A,b_2}^I)$. By linearity of projections we have
\begin{displaymath}
\Pi \cdot \tr_{\mathcal C}(P_{A,b}^I)\subseteq \Pi \cdot
\tr_{\mathcal C}(P_{A,b_1}^I)+
\Pi \cdot \tr_{\mathcal C}(P_{A,b_2}^I).
\end{displaymath}
We will show that $\tr_{\mathcal D}(P_{\bar{A},\bar{b}}^I) \subseteq
\tr_{\mathcal D}(P_{\bar{A},\bar{b}_1}^I) + \tr_{\mathcal
  D}(P_{\bar{A},\bar{b}_2}^I)$. Let $x \in \tr_{\mathcal
  D}(P_{\bar{A},\bar{b}}^I) \subseteq \Pi \cdot \tr_{\mathcal
  C}(P_{A,b}^I))$. Thus there is a representation $x = x_1 + x_2$
with $x_i \in \Pi \cdot \tr_{\mathcal C}(P_{A,b_i}^I)$. We will
show that $x_i \in \tr_{\mathcal D}(P_{\bar{A},\bar{b}_i}^I)$ for
$i=1,2$. Suppose not and let 
w.l.o.g.~$x_1 \notin \tr_{\mathcal D}(P_{\bar{A},\bar{b}_1}^I)$. Then
there is $d \in {\mathcal D}$ with $d \sqsubseteq x_1$. But then $d
\sqsubseteq x_1+x_2 = x$ because $x_i \in \Z^m_+$. This contradicts
the fact that $x \in \tr_{\mathcal D}(P_{\bar{A},\bar{b}}^I)$ and
concludes the proof. \end{proof} 

This implies: If we have a finite generating system for the fibers
$P_{A,b}^I$ (or $\tr_{\mathcal C}(P_{A,b}^I)$) then we have a finite
generating system for the projected fibers $P_{\bar{A},\bar{b}}^I$
(or $\tr_{\mathcal D}(P_{\bar{A},\bar{b}}^I)$) w.r.t.~the monoid
${\mathcal B}$. Therefore we may reformulate the integer Minkowski
program on the projected fibers as an integer linear program, because we
may do so for the original fibers.

\section{The mixed-integer case: Indecomposable Polytopes}
\label{section mixed-integer}
\subsection{The general case}
In this section we will consider finitely generated families of
polyhedral sets depending on integral vectors. We refer to the notation of
\citet{henk-koeppe-weismantel-2000:miptestset, koeppe:diss}. Let $W
\in \Z^{m \times n}$ be a fixed but arbitrary integral matrix with
row vectors $w^i \in \Z^n$, $1 \leq i \leq m$. We assume that 
$$\pos\{w^1, \ldots, w^m\} = \R^n, $$ 
where $\pos$ denotes the positive hull. Thus for every $u \in \R^m$ the
set $P_u = \{y \in \R^n : Wy \leq u\}$ is a polytope. We are
interested in the set of all nonempty polytopes with integral
right-hand side arising in this way. We set:
$${\mathcal U}(W) = \{u \in \R^m : P_u \neq \emptyset \},$$
and consider ${\mathcal U}(W) \cap \Z^m$.

\begin{Definition}
\label{int-decomposability}
A polytope $P_z$, $z \in {\mathcal U}(W) \cap \Z^m$ is called
\emph{integrally decomposable} if there exist $P_{z^1}$, $P_{z^2}$
not homothetic to $P_z$ such that $P_z = P_{z^1} + P_{z^2}$ and $z=
z^1 + z^2$, $z^i \in {\mathcal U}(W) \cap \Z^m$. $P_z$ is called
\emph{integrally decomposable} otherwise.
\end{Definition}

We have the following result in
\citet{henk-koeppe-weismantel-2000:miptestset}: 

\begin{Theorem}
\label{finite-generating-set}
There exist finitely many vectors $h^1, \ldots, h^k \in {\mathcal
  U}(W) \cap \Z^m$ such that for every polytope $P_z$, $z \in
{\mathcal U}(W) \cap \Z^m$, there exist $h^{j_1}, \ldots,
h^{j_{2m-2-n}}$ and non-negative integers $\lambda_{j_1}, 
\ldots, \lambda_{j_{2m-2-n}}$ such that
\begin{displaymath}
P_z = \sum\limits_{i=1}^{2m-2-n}{\lambda_{j_i} P_{h^{j_i}}} \quad
\text{and} \quad z = \sum\limits_{i=1}^{2m-2-n}{\lambda_{j_i} h^{j_i}}.
\end{displaymath} 
\end{Theorem}

This means in particular, that it is possible to model a mixed-integer
Minkowski program on polytopes with integral right-hand side as an
integer linear program according to Theorem \ref{IP und IMP}.

\subsection{Restricted right-hand sides}
In our application to network design the polytopes $P_z$ introduced
in the previous subsection will represent the sets of feasible vectors of
flow of particular networks. Therefore it makes sense to claim
non-negativity of the points in a polytope $P_z$. We assume that our
matrix $W$ is of a special structure: $W = (\bar{W}, - \Id)^T \in
\Z^{(m+n) \times n}$ and we
restrict our attention to the following lattice of right-hand sides:
$$\Psi = \{u = (\bar{u}, 0)^T \in \Z^{m+n} : \bar{u} \in
\Z^m\}.$$
For $u \in \Psi$ we have: 
\begin{displaymath}
P_u = \{y \in \R^n : Wy \leq u\} = \{y \in \R^n_+ :
\bar{W}y \leq \bar{u} \}.
\end{displaymath}
As in the integral case of atomic fibers we are interested in
sublattices of $\Psi$, because of flow-conservation constraints or
non-negativity constraints for demands and capacities for
example. 

\begin{Definition}
Let $\Lambda \subseteq \Psi$ denote a lattice and
$M \subseteq \Psi$ denote a monoid.
\begin{enumerate}[\rm(i)]
\item A polytope $P_u$ with $u \in \Lambda$ is called
  \emph{integrally indecomposable w.r.t.~$\Lambda$} if there are no
  vectors $0 \neq u_1,u_2 \in \Lambda$ with $P_u = P_{u_1} + P_{u_2}$ and $u
  = u_1 + u_2$. 
\item A polytope $P_u$  with $u \in M$ is called
  \emph{integrally indecomposable w.r.t.~$M$} if there are no
  vectors $0 \neq u_1,u_2 \in M$ with $P_u = P_{u_1} + P_{u_2}$ and $u
  = u_1 + u_2$.
\end{enumerate}
\end{Definition}

It was already proved in \citet{henk-koeppe-weismantel-2000:miptestset},
Remark 3.1, that there are only finitely many integrally indecomposable
polytopes w.r.t.~a sublattice $\Lambda$ of $\Z^m$.  We may extend this result
to monoids under certain conditions. Let $M$ be a monoid which is finitely
generated, i.e., which is generated by $m_1, \ldots, m_t \in M$ and let
$\Lambda_M$ the lattice generated by $m_1, \ldots, m_t$. 

\begin{Lemma}
\label{mon-indec-polytopes}
Let $M$ be a monoid and let $W\Z^n_+ \cap \Lambda_M \subseteq M$. Then
there are only finitely many integrally indecomposable polytopes
w.r.t.~$M$. 
\end{Lemma}
\begin{proof}
We know that there are only finitely many integrally indecomposable
polytopes $P_z$ w.r.t.~$\Lambda_M$. As $M \subseteq \Lambda_M$ we
know, that if $z \in M$ and $P_z$ integrally indecomposable
w.r.t.~$\Lambda_M$ then $P_z$ is integrally indecomposable
w.r.t.~$M$. On the other hand, let $P_z = P_{z_1} + P_{z_2}$ with $z
= z_1 + z_2$, $z,z_1,z_2 \in \Lambda_M$ and $P_z, P_{z_i} \neq
\emptyset$. Then $z, z_i \in W\Z^n$, $i=1,2$. As $W\Z^n \cap \Lambda_M =
W\Z^n_+ \cap \Lambda_M \subseteq M$, we have $z,z_1,z_2 \in M$. This
means that the polytopes which are integrally indecomposable
w.r.t.~$\Lambda_M$ are integrally
indecomposable w.r.t.~$M$ and vice-versa. This implies the finiteness
of the number of integrally indecomposable polytopes w.r.t.~$M$. \end{proof} 

\begin{Corollary}
\label{truncation}
Let ${\mathcal C}$ be a set of vectors (not necessarily
finite). Consider the truncated polytopes $\emptyset \neq 
\tr_{\mathcal C}(P_u) = \{y \in \R^n : Wy \leq u, \; \nexists c
\in {\mathcal C} \text{ with } c \sqsubseteq y\}$. There is a finite
generating set for these truncated polytopes.
\end{Corollary}
\begin{proof}
This is a direct consequence of Lemma \ref{trunc subset normal}. \end{proof}

\begin{Corollary}
One may reformulate integer Minkowski programs on (truncated)
polytopes with integral right-hand side as integer linear programs.
\end{Corollary}
\begin{proof} 
This comes from Theorem \ref{tr-reformulation} and Theorem
\ref{IP und IMP}. \end{proof}

\section{Network design problems and atomic fibers}
\label{section formulations}
In this section, we will treat the network design problem in terms
of the analysis of atomic fibers. The outline of this section is as
follows: The first two subsections will introduce two formulations of
the network design problem, the node-arc and the path-cycle
formulation. The following subsection will give the connection between
the path-cycle and the node-arc formulation and introduce the notion
of irreducible networks. The
following two subsections will give the connection between the
irreducible networks and the atomic fibers for both our
formulations.  The last subsection will finally present the network design
problem as an integer Minkowski program and its reformulations as 
an integer linear program.

The reason for studying both formulations is the following.  The path
formulation is easier to study than the node-arc formulation:  In the path
formulation, ordinary atomic fibers appear, whereas in the node-arc
formulation, we need to consider atomic fibers with respect to a monoid of
feasible right-hand side vectors.  On the other hand, the node-arc
formulation is essentially a projection of the path-cycle formulation.
Therefore, the set of atomic fibers is in general 
much smaller in the node-arc formulation.  Consequently, the computation of the atomic fibers is
more efficient, and also the reformulation as an integer linear program is
more compact when we start with the node-arc formulation.
\smallbreak

Let $G=(V,A)$ be the supply digraph, which is connected. Each arc
$a \in A$ of the graph has a capacity $c_a \in \Z_+$ such that $c
\in \Z_+^{|A|}$ and we are given a demand vector
$d=(d_1,\ldots, d_k)$, $d_i \in \Z_+ \; \forall i \in \{1, \ldots,
k\}$. For each commodity $l=1,\ldots, k$ let $s_l \in V$ be its source,
$t_l \in V$ its sink and $\Pi_l$ the set of all paths from $s_l$ to
$t_l$. Such a setting is called a {\em multicommodity network}. We
will denote it by $N=(V,A,d,c)$.

\subsection{Node-arc formulation} \label{netdes}
In this subsection, we will give an exact definition of the
problems we consider in this paper.  Let $f_a^l \in \Z_+$ be the
part of the flow of commodity $l$ which uses arc $a$. The problem of
finding a feasible flow w.r.t.~the capacity vector $c$ and the demand
vector $d$ can be formulated as follows:
\begin{subequations}
\label{node-arc}
\begin{alignat}{4}
\label{na-capacity} \sum\limits_{l=1}^{k}{f_a^l}& & && + s_a &\; = c_a && \qquad \forall a \in A \\
\label{na-demand} \sum\limits_{a \in \delta^+(s_l)}{f_a^l}&- \sum\limits_{a \in
  \delta^-(s_l)}{f_a^l}& && & \; =d_l &&  \qquad \forall l=1, \ldots, k \\
\label{na-flow-conserv} \sum\limits_{a \in \delta^+(x)}{f_a^l}&- \sum\limits_{a \in
  \delta^-(x)}{f_a^l}& && & \; = 0 && \qquad \forall x \in V\setminus \{s_l,t_l\}, \;
\forall l=1, \ldots, k  \\
\label{na-integer}& f_a^l, s_a \in \Z_+  & && & && \qquad \forall a \in A, \, \forall l=1, \ldots, k,
\end{alignat}
\end{subequations}
where $s_a$, $a \in A$, denote the slack variables.
 We denote by $C$ the
matrix corresponding to the system of equations \eqref{node-arc}.

\begin{Definition}
\label{Multicommodity Flow Problem} Let $G=(V,A)$ be a digraph,
$d=(d_1, \ldots, d_k) \in \Z_+^k$ a given demand vectors and $c
\in \Z_+^{|A|}$ a given capacity vector. Let $K_a^l \in \R_+$
denote the costs of routing one unit of commodity $l$ through arc
$a$. The program
\begin{equation}
\label{mcf-na}
 \begin{array}{llcll}
\min & K^Tf \\
\st & C \cdot (f,s)^T & = & (c,d,0)^T\\
& f_a^l,\, s_a \in \Z_+ & & &\forall a \in A, \, \forall l=1, \ldots, k
\end{array}
\end{equation}
is called the {\em multicommodity flow problem} for integer flows
with respect to the node-arc formulation.
\end{Definition}

\noindent In the multicommodity flow problem, we are given a
demand vector and a capacity vector. We want to find a
routing of the demands, which is minimal w.r.t.~the costs of
routing flow. In the network design problem, we are given a vector
of demand and want to find a capacity vector minimizing
installation costs, such that we can find a routing for the
demands.

\begin{Definition}
\label{Network Design Problem} Let $G=(V,A)$ be a digraph,
$d=(d_1, \ldots, d_k) \in \Z_+^k$ a given demand vector. Let
additionally $W_a$, $a \in A$, denote the costs of installing one
unit of capacity on arc $a$. The program
\begin{equation}
\label{nd-na}
 \begin{array}{llcll}
\min & W^T~c & \\
\st & C \cdot (f,s)^T & = & (c,d,0)^T \\
& f_a^l,\, s_a \in \Z_+ & & &\forall a \in A, \, \forall l=1, \ldots, k\\
& c_a \in \Z_+ & & & \forall a \in A
\end{array}
\end{equation}
is called the {\em network design problem} for integer flows with
respect to the node-arc formulation.
\end{Definition}

\subsection{Path-cycle formulation}
\label{path formulation} Let ${\cal W}$ be the set of directed
cycles in $G=(V,A)$. $y_p^l$ denotes the part of commodity $l$
routed on path $p$, $p \in \Pi_l$. $y_w^l$ denotes the part of
commodity $l$, which circulates on cycle $w \in {\cal W}$. The
problem of finding a feasible flow w.r.t.~the capacity vector $c$ and
the demand vector $d$ can be formulated as follows:
\begin{subequations}
\label{pc}
\begin{alignat}{4}
\label{pc-capacity}\sum\limits_{l=1}^{k}{} & \biggl(  \sum\limits_{\substack{p \in
        \Pi_l\\ a \in p}}  {y_p^l} +
 \sum\limits_{\substack{w \in {\cal W}\\a \in w}}{y_w^l} \biggr) & &&
{}+ s_a^l
& \; = c_a && \qquad \forall a \in A \\
\label{pc-demand}& \sum\limits_{p \in \Pi_l} {y_p^l} &  &&  &\; =   d_l && \qquad
\forall l \in \{1, \ldots,k\}\\
\label{pc-integer} &  y_p^l, y_w^l \in \Z_+ & && & && \qquad \forall p \in \Pi_l, \forall l \in \{1,
\ldots, k \}, 
\end{alignat}
\end{subequations}
where $s_a$, $a \in A$ denote the slack variables.

As in the previous subsection we define the network
design problem and the multicommodity flow problem for integer
flows although this time w.r.t.~the path-cycle formulation. With
the notation of the previous subsection, we have the following
formulation of the multicommodity flow problem:
\begin{equation}
\label{mcf-pc}
 \begin{array}{ll}
\min & \sum\limits_{a \in
A}{\sum\limits_{l=1}^{k}{K_a^l\biggl(\sum\limits_{\substack{p \in \Pi_l\\ a \in
p}}{y_p^l} + \sum\limits_{\substack{w \in {\cal W}\\ a \in w}}{y_w^l}\biggr)}}\\
\st & \eqref{pc},\\
\end{array}
\end{equation}
and the following formulation of the network design problem:
\begin{equation}
\label{nd-pc}
 \begin{array}{ll}
\min & W^T~c\\
\st & \eqref{pc} \\
& c_a \in \Z_+ \; \forall a \in A.\\
\end{array}
\end{equation}

\noindent Now we look at the multicommodity flow problem: As we
have non-negative costs associated with the arcs ($K \in \R_+^{|A| \cdot k}$),
there will be no cycle flow in the optimal solution, i.e., $y_w^l =
0$ $ \forall w \in {\cal W}, \, \forall l \in \{1, \ldots, k\}$.
So we can eliminate the cycle-flow variables $y_w^l$. For the
network design problem, too, we can eliminate the cycle-flow
variables $y_w^l$. We obtain the following formulation:
\begin{subequations} \label{p}
\begin{alignat}{3}
 \label{p-demand} \sum\limits_{p \in \Pi_l} {y_p^l}&& &=  d_l && \qquad \forall l \in \{1,
\ldots, k\}\\
 \label{p-capacity} \sum\limits_{l=1}^{k}{}\sum\limits_{\substack{p \in \Pi_l\\ a \in
      p}}  {y_p^l}&& \; + s_a &\leq  c_a && \qquad \forall a \in A \\
\label{p-intger}  y_p^l \in \Z_+ && & && \qquad \forall p \in \Pi_l, \forall l \in \{1, \ldots,
k \}\\
s_a \in \Z_+ && & && \qquad \forall a \in A,
\end{alignat}
\end{subequations}
where $s_a$, $a \in A$, denote the slack variables.

Let $D$ denote the matrix corresponding to the left-hand side of
formulation \eqref{p}. We obtain a new formulation of the
multicommodity flow problem
\begin{equation}
\label{mcf-p}
\begin{array}{lrcll}
\min & \sum\limits_{a \in
A}{\sum\limits_{l=1}^{k}{K_a^l(\sum\limits_{\substack{p \in \Pi_l\\ a \in
p}}{y_p^l})}}\\
\st & D \cdot (y,s)^T & =& (d,c)^T\\
& y_p^l \in \Z_+ & & & \forall p \in \Pi_l, \; \forall l=1, \ldots,
k\\
& s_a \in \Z_+ & & & \forall a \in A
\end{array}
\end{equation}
and a new formulation of the network design problem
\begin{equation}
\label{nd-p}
\begin{array}{lrcll}
\min & W^T~c\\
\st & D \cdot (y,s)^T & =& (d,c) \\
& y_p^l \in \Z_+ & & & \forall p \in \Pi_l \; \forall l=1, \ldots, k\\
& c_a,\, s_a \in \Z_+ & & &\forall a \in A.
\end{array}
\end{equation}
This formulation is called the {\em path formulation}.

\subsection{The node-arc formulation as a projection of the path-cycle formulation}
It is clear that one can convert feasible solutions of the
multicommodity flow problem in the path-cycle formulation to to feasible
solutions of the node-arc formulation via a projection. Let $y$ be a
feasible integer solution of the multicommodity flow problem with respect to
the path-cycle formulation. We set:
$$ f_a ^l := \sum\limits_{\substack{p \in \Pi_l \\ a \in p}}{y_p^l} +
\sum\limits_{\substack{w \in {\cal W} \\ a \in w}}{y_w^l} \qquad \forall a \in
A, \; \forall l \in \{1, \ldots, k\}.$$

Then, $f$ is a feasible integer solution of the multicommodity flow
problem with respect to the node-arc formulation. Indeed, we have
integrality of the components of $f$ because $y$ is integer. The
capacity constraints are respected because:
$$ \sum\limits_{l=1}^{k}{f_a^l} + s_a =
\sum\limits_{l=1}^{k}{\biggl ( \sum\limits_{\substack{p \in \Pi_l \\ a \in p}}{y_p^l}
  + \sum\limits_{\substack{w \in {\cal W} \\ a \in w}}{y_w^l}\biggr )}
+ s_a =
c_a \qquad \forall a \in A.$$
It remains to check whether the flow-conservation constraints are
respected. To this aim, we look at a path $p = (s_l, x_1, \ldots, x_n,
t_l)$  with $y_p^l >0$. The flow on path $p$ respects the
flow-conservation constraints for every $x_i \in p$. This observation
is true for every path $p \in \Pi$.  The same
assertion is valid for all cycles $w$ with $y_w^l >0$. Now we look at
some node $x_i \in V\backslash \{s_l,t_l\}$. We have:
$$\sum\limits_{a \in \delta^+(x_i)}{f_a^l} - \sum\limits_{a \in
  \delta^-(x_i)}{f_a^l} = 0 \qquad \forall l \in \{1, \ldots, k\}. $$
Analogous arguments yield:
$$\sum\limits_{a \in \delta^+(s_l)}{f_a^l} - \sum\limits_{a \in
  \delta^-(s_l)}{f_a^l} = d_l \qquad \forall l \in \{1, \ldots,
k\}. $$
Therefore the path-cycle flow $y$ determines the node-arc flow $f$
uniquely. 

Also, one can convert feasible solutions of the
multicommodity flow problem in node-arc formulation to feasible
solutions in path-cycle formulation.
\begin{Lemma}[see \cite{AhujaMagnantiOrlin93}, Theorem 3.5]
Every non-negative arc flow $f$ can be represented as a path and
cycle flow $y$ (though not necessarily uniquely).
\end{Lemma}
However, the solution is not uniquely determined in general.
\begin{Example} 
  To see an example of this non-uniqueness consider the one-commodity-digraph
  in Figure \ref{Nicht-Eindeutigkeit}. 
\begin{figure}[ht]
  \centering \input{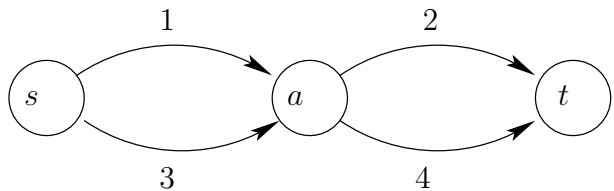} 
  \caption{A network with non-unique path-flow} 
  \label{Nicht-Eindeutigkeit} 
\end{figure} 
There are $4$ paths from the source $s$ to the sink $t$: $p_1 =\{1,2\}$,
$p_2= \{1,4\}$, $p_3=\{3,2\}$ and $p_4=\{3,4\}$. We consider the arc-flow
vector $f=(1,1,1,1)$, i.e., one unit of flow on all arcs. This flow can be
represented in two ways as path-flow: $y^1=(y_{p_1}^1,
y_{p_2}^1,y_{p_3}^1,y_{p_4}^1)=(1,0,0,1)$ or $y^2=(y_{p_1}^2,
y_{p_2}^2,y_{p_3}^2,y_{p_4}^2)=(0,1,1,0)$. 
\end{Example}\smallbreak

\subsection{Irreducible networks in the node-arc and path formulation}
\label{irreducibility in both formulations}
\label{def irreducibility} 
We will now introduce the notion of irreducibility of networks. 
To this aim, we consider the set of feasible solutions of the 
multicommodity flow problem in the node-arc and in the path formulation. 

\begin{Definition}
Let $G=(V,A)$ be a digraph, $d \in \Z^k_+$ the demands and $c \in
\Z^{|A|}_+$ 
the capacities. We denote by ${\cal F}^{\rm{P}}$ the set of feasible integer
solutions of 
the multicommodity flow problem on $N=(V,A,d,c)$ w.r.t.~the path 
formulation. Analogously, we denote by ${\cal F}^{\rm{NA}}$ the set of 
feasible integer solutions w.r.t.~the node-arc formulation.
\end{Definition}

\begin{Definition}
Let $G=(V,A)$ be a digraph, $d \in \Z_+^k$ the demands and $c \in
\Z^{|A|}_+$ 
the capacities. We denote by ${\mathcal FN}^{\rm{P}}$ the set of all non-cyclic
 feasible integer solutions of the multicommodity flow problem on
 $N=(V,A,d,c)$. This means: If $y \in {\mathcal FN}^{\rm{P}}$, then $y$ contains 
no circulation of flow. ${\mathcal FN}^{\rm{NA}}$ is defined analogously.
 Of course we have the following relations:
$${\mathcal FN}^{\rm{P}} \subseteq {\cal F}^{\rm{P}} \qquad {\mathcal FN}^{\rm{NA}} 
\subseteq {\cal F}^{\rm{NA}}.$$
\end{Definition}

To illustrate these definitions, we consider the following
example.
\begin{Example} 
  Let $G=(V,A)$ the digraph in Figure~\ref{circulations}. 

\begin{figure}[ht]
  \centering \input{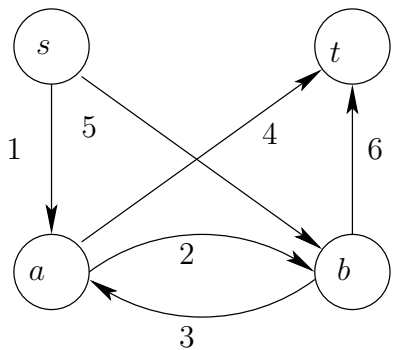} 
  \caption{A digraph with possible circulations} 
  \label{circulations} 
\end{figure} 

We have $4$ paths in this example: $p_1 = \{1,4\}$, $p_2= \{1,2,6\}$,
$p_3=\{5,6\}$ and $p_4=\{5,3,4\}$.  Let $d=2$ and $c=(1,1,1,1,1,1)$. Then $y =
(y_{p_1},y_{p_2},y_{p_3},y_{p_4}) = (0,1,0,1)$ is a feasible integer solution
of the multicommodity flow problem w.r.t.~the path formulation. This means $y
\in {\cal F}^P$.  But $y$ contains a circulation of flow on arcs $2$ and
$3$. Therefore it is not contained in ${\mathcal FN}^{\rm{P}}$. The arc flow
$f$ which is determined by $y$ is contained in ${\cal F}^{\rm{NA}}$ but not in
${\mathcal FN}^{\rm{NA}}$. 
\end{Example} 

Now we can give the definition of irreducibility of networks.

\begin{Definition}
\label{irreducibility}
Let $G=(V,A)$ be a digraph, $d \in \Z_+^k$ the demand and $c \in
\Z^{|A|}_+$ the arc-capacity of a network. A \emph{decomposition
  w.r.t.~the path formulation (w.r.t.~the node-arc formulation)} of
the network $N=(V,A,d,c)$ is given by a decomposition of the
capacity vector and the demand vector $c=c_1+c_2$, $c_1,c_2 \in
\Z_+^{|A|}$, $d= d_1 + d_2$, $d_1,d_2 \in \Z^k_+$, 
 such that the set of non-cyclic
feasible integer solutions of the multicommodity flow problem can be
obtained as the Minkowski-sum: ${\mathcal FN}^{\rm{P}} \subseteq {\cal
  FN}_1^{\rm{P}} +
{\cal FN}_2^{\rm{P}}$ (${\mathcal FN}^{\rm{NA}} = {\mathcal FN}^{\rm{NA}}_1 +
{\mathcal FN}^{\rm{NA}}_2$), where ${\cal FN}_i^{\rm{P}}$ (${\mathcal FN}^{\rm
  {NA}}_i$) is the set of non-cyclic
feasible integer solutions of the multicommodity flow problem on
$N_i=(V,A,d_i,c_i)$, $i=1,2$. If the network can not be decomposed  
in this way, it is called \emph{irreducible} or \emph{indecomposable
  w.r.t.~the path formulation (w.r.t.~the node-arc formulation)}.
\end{Definition}

The notion of
irreducibility depends on the formulation of the multicommodity flow
problem, as the following example illustrates.

\begin{Example} 
  Let us consider the network in Figure \ref{equivalence}, with all the arc
  capacities equal to $1$. 

\begin{figure}[ht]
  \centering \input{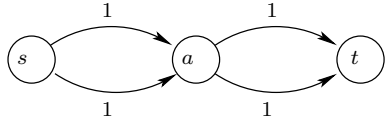} 
  \caption{A network with demand $d = 2$} 
  \label{equivalence} 
\end{figure} 

We have already considered this example in Figure
\ref{Nicht-Eindeutigkeit}. The set of solutions of the multicommodity flow
problem w.r.t.~the node-arc formulation is a singleton: ${\cal F}^{\rm{NA}} =
{\mathcal FN}^{\rm{NA}} = \{(1, 1, 1, 1)\}$. The network is decomposable
w.r.t.~the node-arc formulation. Its possible decompositions are shown in
Figure \ref{equivalence2}. 

\begin{figure}[ht]
  \centering \input{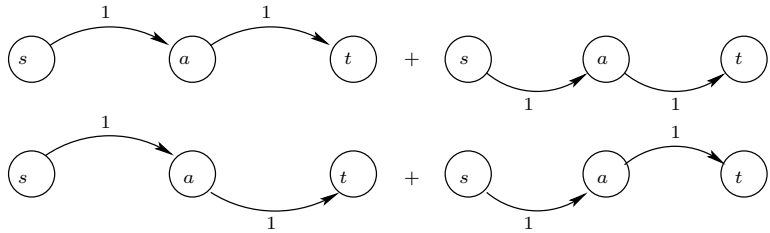} 
  \caption{Decompositions of Figure \ref{equivalence} with demands $d_i =1$,
    $i=1,2$} 
  \label{equivalence2} 
\end{figure} 

In contrast to this, the network in Figure \ref{equivalence} is irreducible
w.r.t.~the path formulation. With the notation of Figure
\ref{Nicht-Eindeutigkeit}, the set of solutions of the multicommodity flow
problem w.r.t.~the path formulation is: ${\cal F}^{\rm {P}} = {\cal
  FN}^{\rm{P}} = \{(1,0,0,1), (0,1,1,0)\}$.  There is no non-trivial
decomposition of the network with the set of solutions of the multicommodity
flow problem decomposing according to it. 
\end{Example} 

Nevertheless we have a connection between irreducibility w.r.t.~the
node-arc formulation and irreducibility w.r.t.~the path
formulation. This connection will be formulated and proved in Theorem
\ref{irred-na-p}.

\subsection{Irreducible networks and atomic fibers in the path formulation}
 Now we want to explore the connection between
irreducible networks and atomic fibers. We recall that a network
is referred to be decomposable w.r.t.~the path formulation for a given vector $d$ of demand and a
given vector $c$ of capacity, if there is a decomposition $d_1,
d_2$ and $c_1, c_2$ of these vectors, such that the feasible non-cyclic
solutions of the multicommodity flow problem decompose according
to it. The network is irreducible w.r.t.~the path formulation otherwise. In fact we would like
to prove that the irreducible networks form a certain subset of the atomic
fibers, i.e., we want to show that they are the truncated fibers
w.r.t.~flow-circulations.  We look at the multicommodity flow problem as defined in
formulation \eqref{mcf-p}. Let $c \in
\Z_+^{|A|}$ be a given capacity vector and $d \in \Z_+^k$ a
given demand vector. 
\begin{equation}
\tag{\ref{mcf-p}}
\begin{array}{lrcll}
\min & \sum\limits_{a \in
A}{\sum\limits_{l=1}^{k}{K_a^l(\sum\limits_{\substack{p \in \Pi_l\\ a \in
p}}{y_p^l})}}\\
\st & D \cdot (y,s)^T & =& (d,c)^T\\
& y_p^l \in \Z_+ & & & \forall p \in \Pi_l, \; \forall l=1, \ldots,
k\\
& s_a \in \Z_+ & & & \forall a \in A
\end{array}
\end{equation}
Let
$$ P_{D,b} := \{(y,s) \in \Z_+^{\sum\nolimits_{l=1}^{k}{|\Pi_l|}}
\times \Z_+^{|A|} : \, D~(y,s)^T = b \}, $$ where $b=(d,c)^T$ is the
right-hand side vector of \eqref{mcf-p}, denote the set of
feasible solutions of the multicommodity flow
problem. Obviously, $P_{D,b}$ is a fiber and the elements $(y,s) \in P_{D,b}$ are
in bijection with the feasible integer
routings of the multicommodity flow problem on $N=(V,A,d,c)$. 
This means: If $y \in {\cal F}^{\rm{P}}$ is a feasible integer solution of
the multicommodity flow problem on $N=(V,A,d,c)$ then $(y,s) \in P_{D,b}$,
where 
$$ s_a := c_a - \sum\limits_{l=1}^{k}{\sum\limits_{\substack{p \in \Pi_l \\ a
      \in p}}{y_p^l}}.$$ 
We have the other direction, too. If $(y,s) \in P_{D,b}$, then $y$ is a
feasible integer routing on $N=(V,A,d,c)$, i.e., $y \in {\cal F}^{\rm{P}}$.

Having this bijection between the set of feasible integer solutions
of the multicommodity flow problem ${\cal F}^{\rm{P}}$ and the elements in
the fiber $P_{D,b}$ we will now establish the connection between
the non-cyclic solutions of the multicommodity flow problem and
truncated fibers. Let $\Gamma = \{\gamma_1, \ldots \gamma_n\}$ denote
the generating set of flow circulations on digraph $G=(V,A)$
w.r.t.~the path formulation. Now, we consider
$\tr_{\Gamma}(P_{D,b})$. 

\begin{Lemma}
\label{irred as trunc fibers}
The elements of the truncated fiber $\tr_{\Gamma}(P_{D,b})$ are in
bijection with the non-cyclic solutions of the multicommodity flow
problem ${\mathcal FN}^{\rm P}$.
\end{Lemma}
\begin{proof}
Suppose $y \in {\mathcal FN}^{\rm P}$. Then $(y,s) \in
\tr_{\Gamma}(P_{D,b})$ with $s$ defined as above. The other direction
is clear as well. \end{proof}

\begin{Corollary}
If $N = (V,A,d,c)$ is an irreducible network w.r.t.~the
path formulation,
 then $P_{D,b}$ with  $b=(c,d)$ is an atomic fiber.
\end{Corollary}
\begin{proof} 
This is a direct consequence of Lemma \ref{trunc subset normal} and
Lemma \ref{irred as trunc fibers}. \end{proof}

\subsection{Irreducible networks and atomic fibers in the node-arc
  formulation}
Up to now we have looked at the atomic fibers of the multicommodity
flow problem with respect to the path formulation. In this subsection
we will treat the atomic fibers of the multicommodity flow problem
with respect to the node-arc formulation.

Consider the the multicommodity flow problem w.r.t.~the node-arc
formulation as defined in \eqref{mcf-na}:
\begin{equation}
\tag{\ref{mcf-na}}
 \begin{array}{llcll}
\min & K^Tf \\
\st & C \cdot (f,s)^T & = & (c,d,0)^T\\
& f_a^l,\, s_a \in \Z_+ & & &\forall a \in A, \, \forall l=1, \ldots, k
\end{array}
\end{equation}
We look at the fibers of the following
form: 
$$P_{C,b} := \{ (f,s) \in \Z^{k\cdot |A|}_+ \times \Z_+^{|A|}
    : C \cdot (f,s)^T
    = b^T\}. $$
If $N=(V,A,d,c)$ is a network, then there is a corresponding non-empty
fiber $P_{C,b}$ with $b=(c,d,0)$. As for the path formulation the
elements in the fiber correspond to the solutions of the
multicommodity flow problem on $N$ and vice-versa. 

But contrarily to the path formulation not all fibers of the matrix $C$
correspond to networks on the digraph $G=(V,A)$. Only a subset of the
fibers corresponds to networks on the digraph $G$. This subset
consists of those fibers with right hand sides $b$ of the form
$b=(c,d,0)$ where $c \in \Z^{|A|}_+$, $d \in \Z_+^k$ and $0$ the
zero-vector. On the other hand all fibers with right hand sides of
this form correspond
to networks of the digraph $G$. This means that we are dealing with a
monoid $M$ of right-hand sides $b$, which can be defined as 
\begin{equation}
\label{themonoid}
M:= \{(c,d,0) \, : \, c \in \Z_+^{|A|}, \, d \in \Z_+^k, \, 0 \in
\Z^{k \cdot (|V|-2)} \}.
\end{equation}

As in the previous subsection, we want to install a connection
between the fibers which are atomic w.r.t.~$M$ and the irreducible
networks. Therefore, let $\Delta = \{\delta_1, \ldots, \delta_m\}$ be
the set of generators of flow-circulations on the digraph $G=(V,A)$
w.r.t.~the node-arc formulation.
Then, we are able to install the connection between irreducible networks
and the fibers which are atomic w.r.t.~$M$. 

\begin{Lemma}
\label{na trunc and irred}
Let $G=(V,A)$ be a digraph, $c \in \Z_+^{|A|}$ the capacity and $d
\in \Z^k_+$ the demand.  $N=(V,A,d,c)$ is an irreducible network if
and only if the corresponding truncated fiber
$\tr_{\Delta}(P_{C,b})$ with $b=(c,d,0)$ is indecomposable  w.r.t.~the
monoid $M$  in \eqref{themonoid}. 
\end{Lemma}
\begin{proof}
The proof is analogous to the proof of Lemma \ref{irred as trunc
  fibers}. \end{proof}

\begin{Corollary}
If $N=(V,A,d,c)$ is an irreducible network w.r.t.~the node-arc
formulation, the $P_{C,b}$ is an atomic fiber w.r.t.~$M$, where
$b=(c,d,0)$. 
\end{Corollary}
\begin{proof} 
This is a consequence of Lemma \ref{trunc subset normal} and of Lemma
\ref{na trunc and irred}. \end{proof}

It remains to show the connection between irreducible networks in the
node-arc and irreducible networks in the path formulation.

\begin{Theorem}
\label{irred-na-p}
A network which is irreducible w.r.t.~the node-arc formulation is     
irreducible w.r.t.~the path formulation.
\end{Theorem}

\begin{proof}
The assertion of the theorem is equivalent to the following
assertion: A network which is decomposable w.r.t.~the path
formulation is decomposable w.r.t.~the node-arc formulation. But this
assertion is a consequence of Lemma \ref{truncated,projected
  fibers}. We have seen in section \ref{irreducibility in both
  formulations} that the node-arc formulation is a projection of the
path-cycle formulation. If we denote $\Lambda = \{\lambda_1, \ldots,
\lambda_t\}$ the generators of flow-circulations w.r.t.~the
path-cycle formulation, then we have: $\Gamma \subseteq \Lambda$ and
we obtain:

A network which is irreducible w.r.t.~the node-arc formulation is
irreducible w.r.t.~the path-cycle formulation. But the non-cyclic
solutions of the multicommodity flow problem w.r.t.~the path-cycle
formulation are in bijection with the non-cyclic solutions of the
multicommodity flow problem w.r.t.~the path formulation. In
particular: The irreducible networks in both formulations are in
bijection and the assertion of the theorem follows. \end{proof}

\subsection{Formulating the network design problem as an integer
  linear program: An example}
\label{IP-Formulation}
 We have seen in the previous subsections that the network design
 problem for integer flows fits in our framework of the integer
 Minkowski programs because there we want to ensure that the set of
 solutions of the multicommodity flow problem for the optimal
 capacity $c$ is not empty. We have also seen that there is a finite
 generating system for the family of feasible solutions sets of
 the multicommodity flow problem. Even though there are no
 Minkowski-additive and antitone functionals involved in the integer
 Minkowski programs corresponding to the network design problem we
 may reformulate it as an integer linear program. As this reformulation
 introduces an integer variable for each irreducible network the
 reformulation process generates an ``extended version'' of the path
 formulation. Note that  a vector $b=(c,d,0)^T$, where $d=e_i$, $i \in \{1,
\ldots,k\}$, and $c= \chi$ the incidence-vector of a path routing
 $d=e_i$ is necessarily contained in the set of irreducible
 networks. Therefore, all path variables $y_p^l$ are contained in the
 reformulation as an integer linear program.

 % Now we can apply the results for atomic fibers to network
%  design. Knowing that there are only finitely many atomic fibers for a
%  given matrix implies the finiteness of the number of irreducible
%  networks for a given digraph. We want to solve the network design
%  problem for integer flows via the irreducible networks.

 We will illustrate the reformulation process by an example.

\begin{Example}
\label{reformulation-example}
 Consider
the digraph $G=(V,A)$ with three arcs and three nodes in Figure~\ref{reformulation example}. 
\begin{figure}[t]
\centering
\input{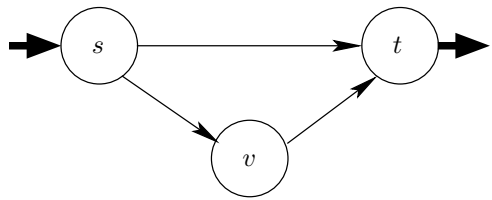}
\caption{The digraph $G=(V,A)$ for one commodity}
\label{reformulation example}
\end{figure} 
For given $d \in \Z_+$ and given $c \in \Z_+^3$ the node-arc
formulation for the digraph in Figure \ref{reformulation example}
looks as follows:
\begin{equation}
\label{na-form}
\begin{array}{c@{\;}c@{\;} c@{\;} c@{\;} c@{\;} c@{\;} c@{\;} c@{\;} c@{\;} c c c c c}
f_1 & + & & & & & & s_1 & & & = & c_1\\
& & f_2 & + & & & &  & s_2 & & = &c_2\\
& & & & f_3 & + & & & &  s_3 & = & c_3 \\
f_1 & + & f_2 & & & & & & &  & = & d\\
& - & f_2 & + & f_3 & & & & & & = & 0
\end{array}
\end{equation}
Let $C$ denote the matrix corresponding to the left-hand side of the
multicommodity flow problem \eqref{na-form}. Then the condition of
finding a feasible integer flow translates to the following condition:
$$P_{C,b} = \{(f,s) \in \Z_+^6 : C \cdot (f,s)^T = b\} \neq
\emptyset,$$
where $b = (c,d,0)^T$. Now the network design problem for the
digraph $G=(V,A)$ for given demand $\bar{d} \in \Z^+$ and vector of
costs $w \in \R_+^3$ may be
formulated as an integer Minkowski program:
\begin{equation}
\label{Minkowski-formulation}
\begin{aligned}
  \mathrm{min}\quad & w^Tc\\
  \mathrm{s.t.}\quad &  d = \bar{d}\\
  & P_{C,b} \neq \emptyset\\
  &  c \in \Z_+^3,  d \in \Z_+
\end{aligned}
\end{equation}

\begin{figure}[ht]
\centering
\input{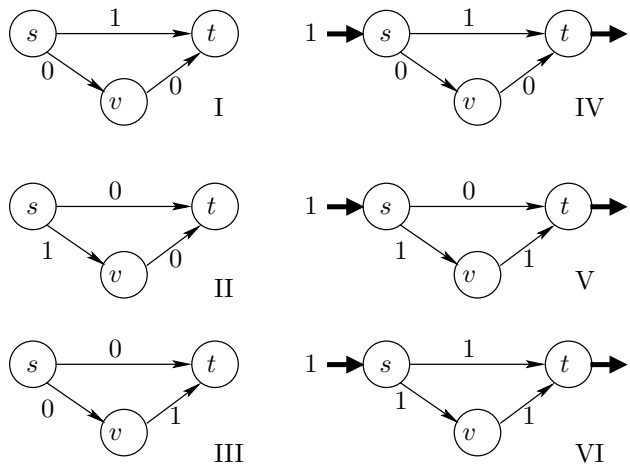}
\caption{The irreducible networks with underlying digraph $G=(V,A)$}
\label{atomic fibers}
\end{figure}

Figure \ref{atomic fibers} shows the fibers of the matrix $C$
which are atomic w.r.t.~the monoid $M=\{(c,d,0) \in \Z_+^6\}$. As
there are no circulations possible on the digraph $G$, these atomic
fibers correspond to the irreducible networks with underlying digraph
$G$.

Now we will reformulate the network design problem with respect to
the digraph $G$ and prescribed demand $\bar{d} \in \Z_+$ as an
integer linear program. Therefore, we introduce a variable $\lambda_i
\in \Z_+$ for all irreducible networks with underlying digraph $G$.
Then the reformulated network design problem is the following one:
\begin{displaymath}
\begin{array}{r@{\;}r@{\,}c@{\,}r@{\,}c@{\,} r@{\,}c@{\,} r@{\,}c@{\,}
    r@{\,} c@{\,}r@{\,}c@{\,} r@{\,}  c}
\rm{min} & w_1\, \lambda_{\rm {I}}& + & w_2\, \lambda_{\rm{II}}& +  & w_3\, \lambda_{\rm{III}}& + &
w_1\, \lambda_{\rm{IV}}& + & (w_2 + w_3) \lambda_{\rm{V}}& + & (w_1 + w_2 + w_3)
\lambda_{\rm{VI}} \\
\rm{s.t.} & & & & & & &\lambda_{\rm{IV}} &+ & \lambda_{\rm{V}} &+& \lambda_{\rm{VI}} & =
\bar{d} \\
& \multicolumn{6}{l}{\lambda_{\rm{I}}, \ldots,  \lambda_{\rm{VI}} \in \Z_+.}
\end{array}
\end{displaymath}
\end{Example}

\section{Design of survivable networks}
\label{section survivability}

An optimal solution of the basic network design problem might give a capacity
installation that admits only one feasible routing.  Clearly, 
if a link or a switching node fails, not all demand can be routed through the
remaining network. In this section,
we will consider various notions of survivability of networks.  
The two basic notions are the arc survivability and the node survivability.
Ensuring arc survivability means to compute a capacity installation such
that there exists a feasible routing of the demand even if an arbitrary 
link of the network fails. 
If a node~$v$ of the network fails, then all arcs incident with
node $v$ fail simultaneously and the commodities incident with node~$v$,
i.e., having sink or source equal to~$v$, fail simultaneously.  Ensuring node
survivability means to compute a capacity installation such
that there exists a feasible routing of the remaining demand even if an
arbitrary node of the network fails.

The two notions can be studied in a unified and generalized failure model.
We refer to the failure model in \cite{Pioro+Medhi:2004}: 
\begin{Definition} 
  A \emph{failure model} is defined as a finite set $\Sigma$ of \emph{failure states}
  of the network.  Each failure state $\sigma\in \Sigma$ is 
  characterized by a vector of \emph{arc-availability coefficients}
  $\alpha_\sigma = (\alpha_{1\sigma}, \ldots, \alpha_{|A|\sigma})$, $0 \leq \alpha_{i,\sigma} \leq
  1$, and a vector of \emph{demand coefficients} $\chi_\sigma = (\chi_{1\sigma} \ldots,
  \chi_{k\sigma})$, $0 \leq \chi_{j\sigma} \leq 1$. The arc-availability coefficients
  $\alpha_{i,\sigma}$ represent the proportion of the capacity of arc $i$ which is
  available in failure state $\sigma \in \Sigma$. The demand coefficients $\chi_{j\sigma}$
  denote the demand of commodity $j$ which must be satisfied in failure state
  $\sigma \in \Sigma$.  
\end{Definition} 
It is obvious that an arc-availability coefficient $\alpha_{i,\sigma}$ is equal to $0$
corresponds to a total failure of arc $i$, whereas $0 < \alpha_{i,\sigma} < 1$
corresponds to a partial failure of this arc and $\alpha_{i,\sigma} =1$ means that
there is no failure of arc $i$. 

In this section we will consider several particular cases for the
arc-availability and the demand coefficients:

\begin{enumerate}
\item $\chi_{\sigma} = \mathbb{1}$ and $\alpha_\sigma \in \{0,1\}^{|A|}$ \newline
This case models a scenario where a single arc or several arcs $a \in
A$ fail and the original demand $d_l$ must be satisfied for all
commodities $l \in \{1, \ldots, k\}$.
\item $\chi_\sigma \in \{0,1\}^k$ and $\alpha_{\sigma} \in \{0,1\}^{|A|}$ with
  $\chi_{i\sigma}=0$ if and only if $\exists v \in V_\sigma \colon v \in \{s_i,t_i\}$ and 
 $\alpha_{a\sigma} =0 $ if and only if $a \in \delta (v)$ for $v \in V_\sigma$, \newline
where $V_\sigma$ denotes the set of nodes failing in failure state $\sigma$.
\item $\chi_\sigma = \mathbb{1}$, $0 \leq \alpha_{i,\sigma} \leq 1$ and
  $\alpha_{j\sigma} \in \{0,1\}$ for $j \neq i$. \newline
This situation corresponds to a partial failure of arc $i$ and either
no failure or total failure of arcs $j \neq i$.
\end{enumerate}

With the help of these ``pure'' scenarios one may also model mixes of
them, e.g., simultaneous failures of several nodes, total failures of
several arcs and a partial failure of a single arc.
Note that within our consideration of arc survivability of networks
we always mean arc survivability w.r.t.~complete rerouting. This
means that there are no routing restrictions imposed in the failure
case: The no-fault and the fault routing may be completely unrelated.

\subsection{Arc survivability w.r.t.~complete rerouting}
In this section we will concentrate on the survivability of networks
w.r.t.~total or partial failures of networks.
The property of a network to be arc-survivable in case of failure of
an arc will be modeled via Minkowski-additive and antitone
functionals depending on the the set of feasible solutions of the
multicommodity flow problem for fixed capacity, demand
and a particular arc of the network. These functionals will be called
\emph{arc-survivability functionals}. Consider the functionals below:
here we define them for the node-arc formulation of the network
design problem. But it is possible to define them in an analogous way
for the path formulation. In the following, the vector $b$ will always
denote $b= (c,d,0)$.  
\begin{equation}
\label{survivability-functions}
 g_a(P_{C,b}) := \min_{(f,s) \in P_{C,b}} \left
   (\sum\limits_{l=1}^{k}{f_a^l} \right ).
\end{equation} 

A network $N=(V,A,d,c)$ is survivable w.r.t.~complete rerouting in case
of total failure of arc $a \in A$ if and only if $g_a(P_{C,b}) =
0$. 

\begin{Lemma}
\label{survivability function}
Let $a \in A$ be fixed.
The functional $g_a(P_{C,b})$ is Minkowski-additive, i.e., if $N=(V,A,d,c)$ is
decomposable into 
$N^1=(V,A,d^1,c^1)$ and $N^2=(V,A,d^2,c^2)$, then we have
$g_a(P_{C,b})=g_a(P_{C,b^1})+g_a(P_{C,b^2})$. Furthermore, $g_a(P_{C,b})$ is an
antitone and non-negative functional.
\end{Lemma}
%%% FIXME: Das kann auch als Folgerung aus Raymonds Atomic-IP-Lemma
%%% betrachtet werden.  Aber das entfernen wir gerade aus dem Fiber-Paper...

\begin{proof} Let $N=(V,A,d,c)$ be a network and $a \in A$ be an arc. Let
$f^*$ be the solution of the multicommodity flow problem
\eqref{mcf-p} on $N$ with minimum aggregated flow on arc $a$. Then we have:
$$g_a(P_{C,b})= \sum\limits_{l=1}^{k}{(f^*)_a^l}.$$ 
W.l.o.g.~we may assume, that $f^* \in {\mathcal FN}^{\rm{NA}}$, i.e., $f^*$ admits
no circulation of flow. 

Now let $N = N_1 + N_2$, implying that $f^* = f_1 + f_2$ with $f_i \in
{\mathcal FN}^{\rm{NA}}_i, \; i=1,2$. This yields 
immediately that the solutions of the multicommodity flow
problem \eqref{mcf-p} on $N_i$, $f_i$, $i=1,2$, admit each minimum
aggregated flow on arc $a$.  Therefore we have:
$$g_a(P_{C,b}) = g_a(P_{C,b^1}) + g_a(P_{C,b^2}).$$ This implies finally the
Minkowski-additivity of the arc-survivability functionals $g_a(P_{C,b})$.  

The arc-survivability functional $g_a(P_{C,b})$ is antitone as the
minimum-function is antitone. Its
non-negativity comes from the non-negativity of feasible flows. \end{proof}

Since $g_a(P_{C,b})$ is Minkowski-additive and antitone we can apply
Theorem \ref{IP und IMP}: We set $\delta := 0$ and add the inequalities
$$ \sum\limits_{i=1}^{t}{\lambda_i~g_a(P_{C,b_i})} \leq 0 \qquad \forall a
\in A$$ to the reformulation of the network design problem as integer
linear program to ensure
survivability w.r.t.~complete rerouting in case of failure of an arc. As the
functionals $g_a(P_{C,b})$ are non-negative we can replace these inequalities
by the following equalities:
$$ \sum\limits_{i=1}^{t}{\lambda_i~g_a(P_{C,b_i})} = 0 \qquad \forall a
\in A$$
and arrive at the following reformulation of the network design
problem under arc survivability constraints as integer linear program:
\begin{equation}
\label{nd-ip-surv}
\begin{array}{llcll}
\min & W^T (\sum\limits_{i=1}^{t}{\lambda_i c_i})\\
\st & \sum\limits_{i=1}^{t}{\lambda_i d_i} & = & d\\
&\sum\limits_{i=1}^{t}{\lambda_i~g_a(P_{C,b_i})} & =&  0 & \forall a
\in A\\
& (c_i,d_i,0)^T = b_i \in F(C)& & & \forall i=1, \ldots, t\\
& \lambda_i \in \Z_+ & & & \forall i = 1, \ldots, t\\
\end{array}
\end{equation}
\begin{Remark} 
  We remark that this integer linear program is infeasible if the topology of
  the network does not support a survivable installation of capacities.  This
  case will occur if there is a commodity pair that does not have two
  arc-disjoint paths. 
\end{Remark} 

We will illustrate this reformulation process with the help of an
example. 
\begin{Example} 
  We resume Example \ref{reformulation-example} and incorporate the
  arc-survivability functionals in our reformulation. Consider Figure
  \ref{reformulation_example:survivability}. The values of the
  arc-survivability functionals on the arcs of the irreducible networks are
  marked in bold. 
\begin{figure}[ht]
  \centering \input{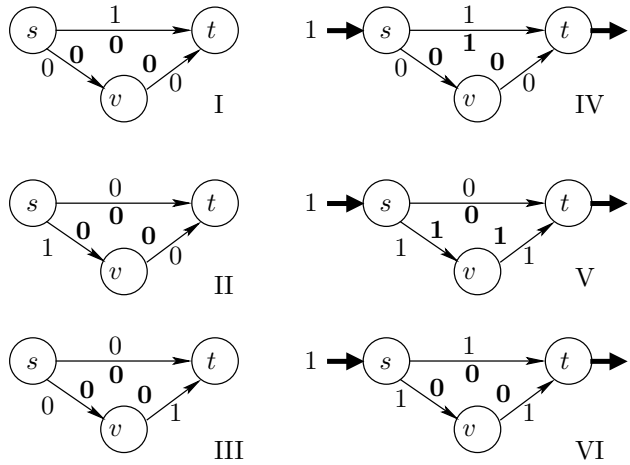} 
  \caption{The values of the arc-survivability functionals} 
  \label{reformulation_example:survivability} 
\end{figure} 

Now we may consider the reformulation of the network design problem under
arc-survivability constraints. For our example this reformulation is the
following one: 
\begin{displaymath} 
  \begin{array}{r@{\;}r@{\,} r@{\,} r@{\,} r@{\,} r@{\,} r@{\,} r@{\,} c}
    \mathrm{min} & w_1\, \lambda_{\mathrm I} + & w_2\, \lambda_{\mathrm{II}} +  & w_3\, \lambda_{\mathrm{III}} + &
    w_1\, \lambda_{\mathrm{IV}} + & (w_2 + w_3) \lambda_{\mathrm{V}} + & (w_1 + w_2 + w_3)
    \lambda_{\mathrm{VI}} \\
    \mathrm{s.t.} & & & & \lambda_{\mathrm{IV}} + & \lambda_{\mathrm{V}}+ & \lambda_{\mathrm{VI}} & =
    \bar{d} \\
    &0\, \lambda_{\mathrm{I}}+ &0\, \lambda_{\mathrm{II}}+ & 0\,
    \lambda_{\mathrm{III}}+ & 1\, \lambda_{\mathrm{IV}}+ & 0\,
    \lambda_{\mathrm{V}}+ & 0\, \lambda_{\mathrm{VI}} & = 0\\
    &0\, \lambda_{\mathrm{I}}+ &0\, \lambda_{\mathrm{II}}+ & 0\,
    \lambda_{\mathrm{III}}+ & 0\, \lambda_{\mathrm{IV}}+ & 1\,
    \lambda_{\mathrm{V}}+ & 0\, \lambda_{\mathrm{VI}} & = 0\\
    &0\, \lambda_{\mathrm{I}}+ &0\, \lambda_{\mathrm{II}}+ & 0\,
    \lambda_{\mathrm{III}}+ & 0\, \lambda_{\mathrm{IV}}+ & 1\,
    \lambda_{\mathrm{V}}+ & 0\, \lambda_{\mathrm{VI}} & = 0\\
    & \multicolumn{6}{l}{\lambda_{\mathrm{I}}, \ldots, \lambda_{\mathrm{VI}} \in \Z_+.}
  \end{array}
\end{displaymath}
\end{Example} 

\begin{Remark}
\label{survivable networks decompose into survivable networks}
Note that if a network is arc-survivable w.r.t.~the failure of an arc
it will decompose into survivable, irreducible summands. This implies
that all variables $\lambda_i$ for non-survivable, irreducible
networks $N_i$ are fixed to $0$ in the reformulation as integer linear
program of the survivable network design problem.
\end{Remark}

\subsection{Partial failures}
Up to now, we have considered total failures of arcs. But it is
possible to model so-called partial failures of arcs. We want to
clarify this with the help of the failure model of
\cite{Pioro+Medhi:2004}. The above discussion has enabled us to model
the $|\Sigma| = |A|$ failure states with $\chi_{\sigma_a} = \mathbb{1}$,
$\alpha_{a\sigma_a}=0$ and $\alpha_{e\sigma_a}=1$ for all $a \neq e \in A$. The
above model guarantees survivability of the network in each of these
failure states. But we may model partial failures of arcs in an
analogous way. To this aim let $\Sigma$ with $|\Sigma| = |A|$ be the index set of the
failure states with $\chi_{\sigma_a} = \mathbb{1}$, $0 < \alpha_{a\sigma_a} < 1$ and
$\alpha_{e\sigma_a}=1$ for all $a \neq e \in A$. Then incorporating the
following inequalities to our model guarantees survivability of the
network in all of these failure states:
\begin{equation}
\label{partial failures}
g_a(c,d) \leq \alpha_{a\sigma_a} c_a.
\end{equation}
These constraints fit in our framework of integer Minkowski programs:
in section \ref{Definitions and first results} we have seen that we
may incorporate non-linear constraints
$$\kappa({\mathcal A}(z)) + h^Tz \leq \delta$$
with $\kappa$ Minkowski-additive and antitone. We already know that
the arc-survivability functionals $g_a(P_{C,b})$ are
Minkowski-additive, antitone and non-negative. Now let $h_a :=
-\alpha_{a,\sigma_a}e_a$ with $e_a \in \Z_+^{|A| + (k-1)|V|}$ the
unit-vector admitting a $1$ in the $a$-th position and let $\tau :=
0$. Then the condition
$$g_a(P_{C,b}) + h^Tb \leq \tau$$
is equivalent to inequality \eqref{partial failures}.

\subsection{Simultaneous failures}
Up to now, we have considered the case of survivability
w.r.t.~complete rerouting in case of failure of only one arc. Our
approach to model survivability of a network in case of failure of an
arc can easily be modified to model survivability w.r.t.~complete
rerouting in case of simultaneous failure of  up to $q$ arcs. In the
model of \cite{Pioro+Medhi:2004} this means: Let $\Sigma$ be an index set
of failure states with $|\Sigma| < \infty$. $\chi_{\sigma} = \mathbb{1}$ for
all $\sigma \in \Sigma$ and $\alpha_{\sigma} \in \{0,1\}^{|A|}$. Now
let $\sigma \in \Sigma$ be 
fixed and let $a_1, \ldots, a_q$ be the set of arcs with
$\alpha_{a_i,\sigma} = 0$. Consider the
Minkowski-additive and antitone functionals presented below:
$$ g_{a_1, \ldots, a_q}(P_{C,b}) := \min_{(f,s) \in P_{C,b}}
\left(\sum\limits_{l=1}^{k} {( f_{a_1}^l + f_{a_2}^l + \ldots +
    f_{a_q}^l )} \right ). $$
We will call functions of the above type \emph{arc-survivability
  functionals for simultaneous failure of arcs $a_1, \ldots, a_q$}.
The Minkowski-additivity of these functionals can be shown analogously
to the proof of Lemma \ref{survivability function}. Antitonicity and
non-negativity of the function follow as well.
As before, a network $N=(V,A,c,d)$ is survivable w.r.t.~complete
rerouting in case of simultaneous failure of
arcs $a_1, \ldots, a_q$ if and only if $ g_{a_1, \ldots, a_q}(P_{C,b}) =
0$. Applying Theorem \ref{IP und IMP} and adding these
equalities to the formulation of the network design problem as
integer Minkowski program for all failure states $\sigma \in \Sigma$ ensures
arc survivability of the network. As
before survivability can be guaranteed by adding linear equalities to
the formulation as integer linear program.

\subsection{Node survivability w.r.t.~complete rerouting}

\begin{figure}[ht]
\centering
\input{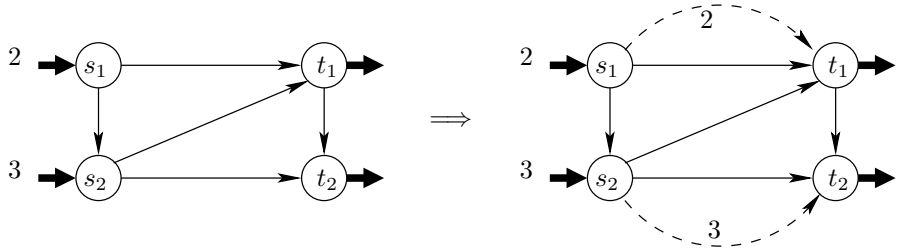}
\caption{Transforming the digraph $G=(V,A)$ to $G'=(V,A \dotcup
  \bar{A})$}
\label{transformation}
\end{figure}

The arc-survivability functionals for simultaneous failure of $q$
arcs will allow us to model node survivability of a network. If a node
$v \in V$ fails then all arcs incident with this node fail
simultaneously. Additionally the demand of all commodities incident
with $v$ fail as well. This implies in particular that the load of the network is
lower in the failure case! 
Therefore we have two different cases: Either node $v$ is the
source or sink of a commodity or it is not. We consider a
transformation of our digraph (see Figure~\ref{transformation}) which
enables us to model node survivability of a network via
arc-survivability functionals.  
For each commodity $l$ we introduce an auxiliary arc $\bar{a}_l =
(s_l, t_l)$. Our digraph now consists of the set of nodes $V$ and the
disjoint union of ordinary arcs $A$ and auxiliary arcs $\bar{A}$,
i.e., $G=(V, A \dotcup \bar{A})$. Each arc $\bar{a}_l \in \bar{A}$ has a
capacity of $c_{\bar{a}_l} = d_l$. 
If a node $v$ fails we want the
commodities incident with this node to be routed via the auxiliary
arcs $\bar{a} \in \delta (v)$.

Now we will model node survivability of a network. To this aim let
$v \in V$ be fixed. We want to model survivability of the network in
case of failure of node $v$. In the model of \cite{Pioro+Medhi:2004}
the corresponding failure state would be the following one:
$\chi_{\sigma} \in \{0,1\}^k$ with $\chi_{\sigma l}= 0$ if and only if $v \in \{s_l,t_l\}$ and
$\alpha_{\sigma} \in \{0,1\}^{|A|}$ with $\alpha_{a\sigma} = 0$ if and only if $a \in \delta
(v)$. Now we want to model survivability of the network in this
failure state. First of all we
have to ensure the existence of a no-fault routing which does not use
any of the auxiliary arcs and satisfies all the demands. This means
that the minimum simultaneous flow on all auxiliary arcs is equal to
$0$: 
\begin{equation}
g_{\bar{a}_1, \ldots, \bar{a}_k}(P_{C,b}) = 0.
\end{equation}

Additionally one has to ensure that there exists a feasible vector of
flow in case of failure of a node $v$. Then the following constraint
ensures survivability of the network in case of failure of arc $v$.
\begin{equation}
g_{\substack{a \in A \cap \delta (v) \\ \bar{a} \in \bar{A}\backslash
    \delta (v)}}(P_{C,b}) = 0
\end{equation}

With the help of these auxiliary arcs and the arc-survivability
functional for simultaneous failure of a set of arcs it is possible
to model node survivability of a network by adding the above
constraints to our formulation of the network design problem as
integer Minkowski program. As shown in the previous section it is
possible to guarantee node survivability by adding linear constraints
to the reformulation as integer linear program.

It is clear that one may model more complex failure situations by
combining and adapting the above scenarios.

%% \section{Computational results}
%% \label{section computations}

\section*{Conclusions}
In a forthcoming paper~\citep*{eisenschmidt-hemmecke-koeppe:network-fibers}, 
we will present algorithms to compute atomic fibers for network design
problems and computational results.
We remark that we expect that the integer linear programs arising from the
reformulation method will be extremely large.  It is therefore desirable to
devise a method for computing atomic fibers on the fly, and to use it within
a specialized branch-cut-and-price algorithm.  However, algorithmic questions like this
are beyond the scope of this paper; they will be the topic of future research.

%% \clearpage
%% \enlargethispage{4ex}
\bibliographystyle{plainnat}
\bibliography{iba-bib,francetelecom}
\end{document}